\documentclass[reqno,10pt]{amsart}
\usepackage{amsmath,amssymb,mathrsfs,amsthm,amsfonts}
\usepackage[usenames,dvipsnames]{xcolor}
\usepackage{comment}
\usepackage{hyperref}
\usepackage{graphicx}
\usepackage{enumerate}
\usepackage{stmaryrd}
\usepackage[all]{xy}
\usepackage{bbm}
\usepackage{hyperref}
\usepackage[capitalise]{cleveref}
\usepackage{caption}
\usepackage{subcaption}
\hypersetup{
     colorlinks=true,
     linkcolor=blue,
     filecolor=black,
     citecolor = blue,      
     urlcolor=blue,
     }
\usepackage{fullpage}

\usepackage{tikz}
\usetikzlibrary{decorations.markings}
\usetikzlibrary{arrows}
\usetikzlibrary{calc}

\newtheorem{theorem}{Theorem}[section]
\newtheorem{lemma}[theorem]{Lemma}

\theoremstyle{definition}

\newtheorem{assumption}[theorem]{Assumption}

\numberwithin{equation}{section}

\graphicspath{{Images/}}

\theoremstyle{plain}
\newtheorem{thm}{Theorem}

\newtheorem{conj}[thm]{Conjecture}
\newtheorem{prop}[thm]{Proposition}

\theoremstyle{definition}
\newtheorem{defn}[thm]{Definition}

\theoremstyle{remark}

\numberwithin{equation}{section}
\numberwithin{thm}{section}

\begin{document}
\title{Convergence of the open WASEP stationary measure without Liggett's condition}
\author[Himwich]{Zoe Himwich}
\address{Department of Mathematics, Columbia University, 2990 Broadway, New York, NY 10027, USA}
\email{himwich@math.columbia.edu}
\date{\today}
\begin{abstract}
    We demonstrate that Liggett's condition can be relaxed without disrupting the convergence of open ASEP stationary measures to the open KPZ stationary measure. This is equivalent to demonstrating that, under weak asymmetry scaling and appropriate scaling of time and space, the four-parameter Askey-Wilson process converges to a two-parameter continuous dual Hahn process. We conjecture that the convergence of the open ASEP height function process to solutions to the open KPZ equation will hold for a wider range of ASEP parameters than those permitted by Liggett's condition. 
\end{abstract}
\maketitle
\setcounter{tocdepth}{1}{
  \hypersetup{linkcolor=black}
  \tableofcontents
}

\section{Introduction}
\subsection{Preface} The open asymmetric simple exclusion process (ASEP), an interacting particle system of great interest in probability and mathematical physics, first appeared in biology, in a 1968 work of MacDonald, Gibbs, and Pipkin~\cite{MGP68} studying the kinetics of protein synthesis. In subsequent years, it also appeared in an early study of exclusion processes by Spitzer~\cite{S70} and Liggett~\cite{Lig75}. The open ASEP is a Markov process which models particles moving on a one-dimensional lattice of finite size with open boundaries, allowing particles to enter and exit the lattice on either side. The process takes its name from its properties: it is \textit{open} in the sense that particles may enter and exit at either boundary, the left and right rates of particle movement are \textit{asymmetric}, meaning that particles move to the right at a higher rate than that at which they move to the left; and only one particle may occupy any location at a given time, a property called \textit{exclusion}. There are other, well studied variants of this system, including a similar model on an infinite lattice, a half infinite lattice, or a periodic domain. 

The open ASEP is a relatively simple model of a system with nonzero net current in its equilibrium state, and, despite its simple definition, it exhibits complex behaviors associated to many-particle systems, such as phase transitions. Due to these properties, the open ASEP has been of great interest to researchers in probability and mathematical physics. Bertini and Giacomin~\cite{BG97} discovered that the ASEP on the full line converges to the Kardar-Parisi-Zhang (KPZ) equation on the full line under weak asymmetry scaling assumptions (see \cref{s:intro2} for a longer discussion of this relationship). Corwin and Shen~\cite{CS18} then proved that the open ASEP on two additional domains (half-open and open ASEP) converge to the half-space and open KPZ equations, respectively, under triple point scaling of the boundary parameters and weak asymmetry scaling. Parekh~\cite{P19} extended the results from \cite{CS18} which were restricted to a certain range of parameters of ASEP to hold for all parameters. As discovered in \cite{BG97}, the open ASEP is closely related to the Kardar-Parisi-Zhang (KPZ) equation (see \cref{s:intro2} for a longer discussion of this relationship). Corwin and Knizel \cite{CK21} used the connection between the open ASEP and the open KPZ equation to give the first construction of the stationary measure for the open KPZ equation through a scaling limit of open ASEP stationary measures. Their work built on extensive literature about the open ASEP, developed over the preceding decades. Work by Derrida, Evans, Hakim, and Pasquier \cite{DEHP93} introduced an algebraic technique, called the matrix product ansatz (MPA) for obtaining the stationary measure of the open ASEP. Uchiyama, Sasamoto, and Wadati \cite{USW} used the MPA expression for the stationary measure of open ASEP to express the measure in terms of Askey-Wilson polynomials. More recently, Bryc and Wesolowski \cite{BWe17} demonstrated a formula for the multi-point correlation function of the stationary measure in terms of a process called the Askey-Wilson Process (which has marginal and transitional laws determined by Askey-Wilson polynomials). \cite{CK21} characterized the Laplace transform of the finite dimensional marginals of the stationary distribution of the open KPZ equation in terms of a stochastic process called the continuous dual Hahn process (see \cref{ss:CDH}), which is the limit of the Askey-Wilson process under the scaling assumptions which take the open ASEP stationary measure to the stationary measure of the open KPZ equation (see \cref{d:scale}).

The open KPZ equation has a unique stationary measure, which was proved by Knizel and Matetski~\cite{KM22} for $u+v>0$ and $\min(u,v)>-1$, and by Parekh~\cite{P22} for all $u,v\in\mathbb{R}$. \cite{CK21} characterized this stationary measure in terms of its Laplace transform for all $u+v\geq 0$. Barraquand and Le Doussal~\cite{BLD22} and Bryc, Kuznetsov, Wang, and Wesolowski~\cite{BKWW23} showed that the stationary measure of the open KPZ equation constructed in \cite{CK21} has a probabilistic description as a process $\{H_{u,v}(x)\}_{x\in[0,1]}=\{B(\frac{x}{2}) - X_{u,v}(x)\}_{x\in[0,1]}$, where $B(x)$ is standard Brownian motion and the law of $X_{u,v}(x)$ is an exponential re-weighting of the law of a Brownian motion by a Radon-Nikodym derivative \begin{align*}
    \frac{d\mathbb{P}_{X}}{d\mathbb{P}_{B}}(f) & = \frac{1}{\mathcal{Z}_{u,v}}e^{-2vf(1)}\left(\int_{0}^{1}e^{-2f(t)}dt\right)^{-(u+v)}.
\end{align*}
Where $f\in C[0,1]$. For the purpose of this paper, we will use the description of the open KPZ equation in terms of its multi-point Laplace transform, see \eqref{e:laplacecharacterization}. 
Recently, work by Barraquand, Corwin, and Yang~\cite{BCY23} constructed the KPZ stationary measure for all $u,v\in\mathbb{R}$ as a scaling limit from the log-gamma directed polymer model, modulo a result extending the existing understanding of convergence of the spacetime process of the log-gamma polymer to convergence on an interval $[0,L]$ (see Claim 4.2~\cite{BCY23}).

The construction \cite{CK21} of the open KPZ stationary measure assumed a relation, called Liggett's condition, on the ASEP parameters. This relation was also used in earlier works~\cite{CS18,P19} that showed that the open ASEP height function process converges as a space-time process to solutions of the open KPZ equation. The contribution of this paper is to demonstrate that Liggett's condition is not necessary to obtain convergence of the stationary measures of open ASEP to that of open KPZ equation. Based on this result, we conjecture that the convergence of open ASEP to open KPZ equation as a space-time process will also hold without this assumption. 

This paper also demonstrates an interesting property of the Askey-Wilson process and its limiting process, the continuous dual Hahn process. The Askey-Wilson process is defined in general with parameters $(A,B,C,D,q)$ (\cref{d:awprocess}), but Liggett's condition fixes the values of $B$ and $D$ as functions of $q$. It is an interesting observation about the structure of the Askey-Wilson and continuous dual Hahn processes, which have been studied in works such as \cite{BW10,BWW23,BW19,B22,WWY23}, to note that the parameters $B$ and $D$ have no effect on the limit process. This observation was anticipated by Bryc, Wang, and Wesolowski~\cite{BWW23} (see, for example, the remarks before Assumption 1.1 in this paper). 

\subsection{The Open ASEP}
To better understand Liggett's condition, we begin by giving a more detailed introduction of the open ASEP.
\begin{figure}[h]
    \centering
    \begin{tikzpicture} \draw (0,0) -- (10,0);
    \draw (-0.5,0) circle (0.5cm);
    \draw (-0.5,-0.75)  node[anchor=north east] {Left Reservoir};
    \draw (0.5,-0.75)  node[anchor=north] {$\alpha$};
    \draw (0.5,0.75)  node[anchor=south] {$\gamma$};
    \draw (9.5,-0.75)  node[anchor=north] {$\delta$};
    \draw (9.5,0.75)  node[anchor=south] {$\beta$};
    \draw (4.5,0.5)  node[anchor=south] {$q$};
    \draw (5.5,0.5)  node[anchor=south] {$1$};
    \draw (10.5,-0.75)  node[anchor=north west] {Right Reservoir};
    \draw[->] (5,0) arc (0:170:0.5cm);
    \draw[->] (9,0) arc (180:30:.60cm);
    \draw[->] (10.1,-0.3) arc (330:190:.60cm);
    \draw[->] (1,0) arc (0:150:.60cm);
    \draw[->] (-.1,-0.3) arc (210:350:.60cm);
    \draw[->] (5,0) arc (180:10:0.5cm);
    \draw (10.5,0) circle (0.5cm);
    \filldraw[fill=black!, draw=black] (1,0) circle (0.1cm);
    \filldraw[fill=white!, draw=black] (2,0) circle (0.1cm);
    \filldraw[fill=black!,draw=black] (3,0) circle (0.1cm);
    \filldraw[fill=white!, draw=black] (4,0) circle (0.1cm);
    \filldraw[fill=black!,draw=black] (5,0) circle (0.1cm);
    \filldraw[fill=white!, draw=black] (6,0) circle (0.1cm);
    \filldraw[fill=white!, draw=black] (7,0) circle (0.1cm);
    \filldraw[fill=white!, draw=black] (8,0) circle (0.1cm);
    \filldraw[fill=black!,draw=black] (9,0) circle (0.1cm);
    \end{tikzpicture}
    \caption{The open ASEP lattice with right, left, and boundary jump rates}
    \label{f:1}
\end{figure}
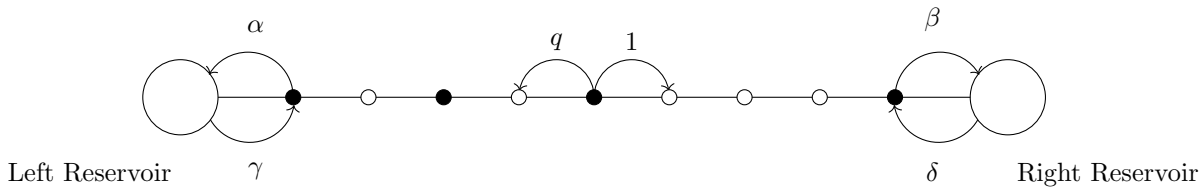

The underlying lattice of the open ASEP has finite size $N$. We will now define the Markov process of particle movement on this lattice in terms of transition rates.
We use $\alpha,\gamma$ to denote the rate of particle movement into and out of the leftmost lattice point from, respectively to, the left reservoir. Similarly, we use $\beta,\delta$ to denote the rate of particle movement out of, respectively into, the rightmost lattice point of the lattice to, respectively from, the right reservoir. Within the lattice, particles jump at times which are determined by independent, exponentially distributed random variables, with rate  $1$ for movement to the right and rate $q\in(0,1)$ for movement to the left. The state space of the open ASEP is $\{0,1\}^{N}$. We refer to elements of the state space by $\tau = (\tau_{1},...,\tau_{N})$ and use the notation $\tau_{i}(t)$ to denote the occupancy of lattice point $i$ at time $t$, with $\tau_{i}(t)=1$ if and only if there is a particle in position $i$ at time $t$. The open ASEP process $\tau(t)$ is a Markov process on the state space $\{0,1\}^{N}$ and the infinitesimal generator $\mathcal{L}(\cdot)$ which acts on functions of the open ASEP state space $f:\{0,1\}^{N}\to\mathbb{R}$ by
\begin{align*}
    (\mathcal{L}f)(\tau) & = \left(\alpha(1-\tau_{1})+\gamma\tau_{1}\right)\cdot\left(f(\tau^{(1)})-f(\tau)\right) \\ & + \sum_{i=1}^{N-1} \left(\tau_{i}(1-\tau_{i+1})+q(1-\tau_{i})\tau_{i+1}\right)\cdot\left(f(\tau^{(i,i+1)})-f(\tau)\right) \\ & + \left(\beta\tau_{N} + \delta(1-\tau_{N})\right)\cdot\left(f(\tau^{(N)})-f(\tau)\right).
\end{align*} In this equation, the notation $\tau^{(i)}$ means that, in $\tau$, the value of $\tau_{i}$ is sent to $1-\tau_{i}$ and the notation $\tau^{(i,i+1)}$ means that $\tau_{i}$ and $\tau_{i+1}$ are swapped. We will denote the stationary measure of the open ASEP with parameters $\alpha,\beta,\gamma,\delta,q$ and lattice size $N$ as $\pi_{N}(\cdot)$, suppressing the dependence on the boundary parameters and the asymmetric jump rate. For the duration of this work, we will use $\langle f\rangle_{N}$ to denote the expectation of a function $f$ under $\pi_{N}(\cdot)$. The time evolution of the open ASEP induces a Markov process, called the height function process, $$h_{N;\alpha,\beta,\gamma,\delta}(t,x) = h_{N;\alpha,\beta,\gamma,\delta}(t,0) + \sum_{i=1}^{\lfloor Nx\rfloor}(2\tau_{i}(t)-1),$$ where $h_{N;\alpha,\beta,\gamma,\delta}(t,0)$ is the net number of particles which have entered through the left boundary at time $t$. If $\tau$ is distributed according to $\pi_{N}(\cdot)$, we refer to the law of the process $\{h_{N;\alpha,\beta,\gamma,\delta}(t,x)-h_{N;\alpha,\beta,\gamma,\delta}(t,0)\}_{x\in [0,1]}$ as the stationary measure for the increment process of the open ASEP height function. In this case the law is independent of $t$, and we will also denote the process by $\{h_{N;\alpha,\beta,\gamma,\delta}(x)-h_{N;\alpha,\beta,\gamma,\delta}(0)\}_{x\in[0,1]}$.

One other standard quantity associated to ASEP is the stationary current, or the average net current of particles entering the ASEP system. As the size of the lattice grows, the value of the limiting stationary current depends on two parameters, $\rho_{\ell}$ and $\rho_{r}$, which are effective densities at the left and right boundaries. Changing parameters to $(A,B,C,D,q)$,
\begin{align*}
     A := \kappa^{+}(q,\beta,\delta), & & B := \kappa^{-}(q,\beta,\delta), & & C := \kappa^{+}(q,\alpha,\gamma), & & D := \kappa^{-}(q,\alpha,\gamma),
 \end{align*}
 \begin{align*}
     \kappa^{\pm}(q,x,y) = \frac{1}{2x}\left(1-q-x+y\pm \sqrt{(1-q-x+y)^{2}+4xy}\right),
 \end{align*}
 we can express $\rho_{\ell}$ and $\rho_{r}$ as
\begin{align*}
    \rho_{\ell} = \frac{1}{1+C}, & & \rho_{r} = \frac{A}{1+A}.
\end{align*} We note that $\rho_{\ell}$ and $\rho_{r}$ are independent of the parameters $B$ and $D$. For a thorough description of the dynamics of open ASEP in terms of parameters $\rho_{\ell}$ and $\rho_{r}$, see the survey \cite{C22}. 

The open ASEP stationary measure is known to converge to the open KPZ stationary measure as $N\to\infty$ and under the scaling
\begin{align}\label{e:scaleold}
    A = q^{v}, & & B = -q, & & C= q^{u}, & & D= -q, & & q=\exp{\left(\frac{-2}{\sqrt{N}}\right)}, 
\end{align} and a particular choice of scaling for the height function (\cref{d:scale}). This choice of scaling for $q$ is called weak asymmetry scaling. \cite{CK21} showed that for any $u,v\in\mathbb{R}$, the $N$-indexed collection of open ASEP stationary measures, scaled according to the weak asymmetry scaling in $N$, was tight, and that sub-sequential limits were stationary measures of the open KPZ equation. Under the constraint that $u+v>0$, they also showed that all sub-sequential limits were equal, and gave a formula for the stationary measures in terms of their Laplace transforms. Later results~\cite{KM22,P22} demonstrated that the open KPZ stationary measure is unique.

The scaling of $A,B,C,D$ with respect to $q$ was chosen to preserve two properties. It preserves a $N^{-\frac{1}{2}}$ scaling around the triple point (discussed in \cite{C22}), 
\begin{align*}
    \rho_{\ell} = \frac{1}{2}+ \frac{u}{2}N^{-\frac{1}{2}} + o(N^{-\frac{1}{2}}), & & \rho_{r} = \frac{1}{2} - \frac{v}{2}N^{-\frac{1}{2}} + o(N^{-\frac{1}{2}}),
\end{align*} by requiring $A=q^{v}$ and $C=q^{u}$ for some $u,v\in\mathbb{R}$. It also preserves a property called Liggett's condition. Liggett's condition is a constraint on the boundary parameters of the system, which, in the parameterization $(A,B,C,D)$, implies that $B=D=-q$. It is typically stated in terms of $(\alpha,\beta,\gamma,\delta)$ as
\begin{align}\label{e:liggett} \alpha + \gamma/q = 1, & & \beta + \delta/q = 1.
\end{align}  

The contribution of this paper is to observe that, without assuming Liggett's condition, the scaling limit of the stationary measure of the open ASEP still converges to the stationary measure of the open KPZ equation. The new scaling assumptions are collected below. 
\begin{assumption}\label{d:scalenew} Our scaling assumptions are
\begin{enumerate}
    \item Weak asymmetry scaling: $q=\exp{\left(\frac{-2}{\sqrt{N}}\right)}.$
    \item Scaling of boundary parameters \begin{align}\label{e:params}
    A  = q^{v}, & & B = -q^{w}, & & C= q^{u}, & & D= -q^{r},
    \end{align} where $u+v>0$ and $w,r>0$.
   \item Height function scaling: For $h_{N}(\cdot)$ the height function of the open ASEP, and $x\in[0,1]$ we define \begin{align*}
       h^{(N)}_{u,v,w,r}(x) & := N^{-\frac{1}{2}}h_{N;\alpha,\beta,\gamma,\delta}\left(Nx\right) .
   \end{align*} We will often drop the subscript and simply use the notation $h^{(N)}(x).$
\end{enumerate}
\end{assumption}

Condition $(2)$ is interchangeable with the following condition on boundary rates $(\alpha,\beta,\gamma,\delta)$,
\begin{equation}\label{e:alsoparams}
\begin{aligned}
    \alpha = \frac{1-q}{(1-q^{r})(1+q^{u})}, & & \beta = \frac{1-q}{(1-q^{w})(1+q^{v})},\\ \gamma = \frac{q^{u+r}(1-q)}{(1-q^{r})(1+q^{u})}, & &  \delta = \frac{q^{v+w}(1-q)}{(1-q^{w})(1+q^{v})}. 
\end{aligned} 
\end{equation} In the equation above, imposing Liggett's condition would set $w=r=1$ and result in a factor of $(1-q)$ canceling from the numerator and denominator in each expression. 
\subsection{Statement and Discussion of the Theorem}\label{s:intro2}

The open KPZ equation is a stochastic partial differential equation defined for $t\in\mathbb{R}_{\geq 0}$ and $x\in[0,1]$ by $$\partial_{t} H(t,x) = \frac{1}{2}\partial^{2}_{x}H(t,x) + \frac{1}{2}\left(\partial_{x}H(t,x)\right)^{2}+\xi(t,x),$$ where $\xi(t,x)$ is space-time white noise. Solutions of the open KPZ equation satisfy \textit{Neumann boundary conditions} $(\partial_{x}H)(t,0) = u$ and $(\partial_{x}H)(t,1) = -v$. This definition is slightly informal. For background information on the open KPZ equation and a careful definition, we direct the reader to \cite{CS18,P19}. Both of those works proved convergence of the open ASEP height function process to Hopf-Cole solutions to the open KPZ equation for a range of boundary parameters under the scaling conditions in \cref{d:scale}, below.

The unique open KPZ stationary measure is the law of the process $\{H_{0}(x)-H_{0}(0)\}_{x\in[0,1]}$ on the space $(D[0,1],\mathcal{B}(D[0,1]))$ such that the law of the space-time process $\{H(t,x)-H(t,0)\}_{t\in\mathbb{R}_{\geq 0},x\in[0,1]}$, produced by evolving initial data $H(0,x)-H(0,0)$ (with $H_{0}(x)=H(0,x)$) according to the open KPZ equation, is invariant in time. \cite{CK21} gave an explicit formula (for $u+v\geq0$) for the open KPZ stationary measure in terms of the Laplace transform of the continuous dual Hahn process \eqref{e:laplacecharacterization}. We will denote the stationary measure of the open KPZ equation with boundary parameters $u,v\in\mathbb{R}$ as the law of the process $\{H_{u,v}(x)-H_{u,v}(0)\}_{x\in[0,1]}$. Using this notation, we can state the main result of this article.

\begin{thm}\label{t:1}
    Assume that the open ASEP parameters and height function satisfy the scaling parameters in \cref{d:scalenew}. The scaling induces a tight sequence of laws of $\{h_{u,v,w,r}^{(N)}(\cdot)-h_{u,v,w,r}^{(N)}(0)\}_{N\in\mathbb{N}}$ on $(C[0,1],\mathcal{B}(C[0,1]))$ and for all $u+v>0$ and $w,r>0$, this sequence converges weakly to the law of $H_{u,v}(\cdot)-H_{u,v}(0)$, the stationary measure of the open KPZ equation, a limit which does not depend on $w,r$. 
\end{thm}

It is natural to consider which of these scaling assumptions is necessary to obtain the convergence of the open ASEP space-time process to that of the open KPZ equation. Convergence results for various ASEP models to KPZ variants~\cite{CS18,P19,Y20} have followed the framework of the foundational paper \cite{BG97}, which approached the problem via the Hopf-Cole transform. Consequently, the argument in each of these works requires Liggett's condition. While Gon\c calves, Perkowski, and Simon~\cite{GPS20} do not use the Hopf-Cole transform, they still assume Liggett's condition. It is not known whether convergence as processes occurs for any choice of $u,v,w,r$ satisfying $u+v>0$ and $w,r>0$. The effective boundary densities, $\rho_{\ell}$ and $\rho_{r}$, depend only on $A$ and $C$, which suggests we may vary $B$ and $D$ without disrupting the triple point scaling. Our result in \cref{t:1} provides further evidence that it is possible to obtain convergence as space-time processes without Liggett's condition, since it is possible at the level of stationary measures. We conjecture that convergence as space-time processes holds using the scaling assumptions in \cref{d:scalenew}.
\begin{assumption}\label{d:scale} The scaling assumptions for the convergence of the open ASEP height function $h_{N}(\cdot)$ to the Hopf-Cole solution to the open KPZ equation are given by condition $(1)$ (weak asymmetry scaling) from \cref{d:scalenew}, with the added conditions 
\begin{enumerate}
    \item Liggett's condition:
   \begin{align*}
       \alpha + \gamma/q=1, & & \beta + \delta/q=1.
   \end{align*}
   \item 4:2:1 Height function scaling: For $h_{N}(\cdot)$ the height function process of the open ASEP, $t\geq 0$ and $x\in[0,1]$ we define \begin{align*}
       h^{(N)}_{u,v}(t,x) & :=N^{-\frac{1}{2}}h_{N;\alpha,\beta,\gamma,\delta}\left(\frac{1}{2}e^{N^{-\frac{1}{2}}}N^{2}t,Nx\right) + \left(\frac{1}{2}N^{-1}+\frac{1}{24}\right)t, \\
       z^{(N)}_{u,v}(t,x) & := \exp{\left(h_{u,v}^{(N)}(t,x)\right)}.
   \end{align*}
   \item H\"older bounds on initial data: The open ASEP initial data is $N$-indexed and satisfies H\"older bounds. For all $\theta\in\left(0,\frac{1}{2}\right)$ and $n\in\mathbb{Z}_{n\geq1}$, there exist constants $C(n)$ and $C(\theta,n)>0$ such that for all $x,x'\in[0,1]$ and all $N\in\mathbb{Z}_{\geq 1}$, \begin{align*}
       \mathbb{E}\left[|z^{(N)}_{u,v}(x)|^{n}\right]^{1/n}\leq C(n), & & \mathbb{E}\left[|z^{(N)}_{u,v}(x)-z^{(N)}_{u,v}(x')|^{n}\right]^{1/n}\leq C(\theta,n)|x-x'|^{\theta},
   \end{align*} where $\mathbb{E}[\cdot]$ is the expectation over ASEP height function initial data $h_{N}(x)$.
\end{enumerate}
\end{assumption}

\begin{conj}\label{c} We denote the Hopf-Cole solution to the open KPZ equation with Neumann boundary parameters $u$ and $v$ by $H_{u,v}(t,x)$. We define $Z_{u,v}(t,x):=\exp{\left(H_{u,v}(t,x)\right)}$. We denote the diffusively scaled height function process with boundary rates $\alpha,\beta,\gamma,\delta$ given by \eqref{e:alsoparams} as $h_{u,v,w,r}^{(N)}(t,x)$, and $z_{u,v,w,r}^{(N)}(t,x):=\exp{\left(h_{u,v,w,r}^{(N)}(t,x)\right)}$. Under weak asymmetry scaling assumptions (\cref{d:scalenew}), and using the height function process scaling defined in \cref{d:scale} (4) and requiring the same H\"older bounds described in \cref{d:scale} (5), then, if the initial data $z^{(N)}_{u,v,w,r}(0,x)$ converges to $Z_{u,v}(0,x)$, $z^{(N)}_{u,v,w,r}(t,x)\implies Z_{u,v}(t,x)$ as a space-time process.
\end{conj}

\subsection{Approach to proving \cref{t:1}}
The proof of the tightness result in \cref{t:1} is an application of an attractive coupling of ASEP proved in \cite{CK21}. This result allows us to obtain H\"older bounds on $h^{(N)}_{u,v,w,r}(x)-h^{(N)}_{u,v,w,r}(0)$ which consequently allow us to conclude that the sequence is tight. For the details of this approach, see \cref{s:proof}.

In order to prove that the sequence of measures in \cref{t:1} converges to the law of $H_{u,v}(\cdot)-H_{u,v}(0),$ we rely on the fact that the open ASEP stationary measure can be written in terms of a stochastic process called the Askey-Wilson process. Consequently, we can describe the Laplace transform of the height function process in terms of the Askey-Wilson process and analyze the asymptotics of that expression in terms of Askey-Wilson densities. We also rely on the fact that pointwise convergence of Laplace transforms on an open set implies convergence in distribution (see \cref{t:laplace}). The main technical contribution of this paper is to show that in a more general range of parameters, for any $w,r>0$ and $u+v>0$, the Askey-Wilson process tangent process in parameters $v,w,u,r$ (see \cref{s:atoms}) converges under \cref{d:scalenew} to the continuous dual Hahn process with parameters $u,v$. Consequently, when $u+v>0$ and $w,r>0$, the tight sequence of laws from \cref{t:1} has a unique limit, which is the stationary measure of the open KPZ equation.

\subsection{Organization of the Paper}
In \cref{s:def}, we introduce the Askey-Wilson process and the continuous dual Hahn process, and state results and definitions pertaining to the asymptotics of $q$-Pochhammer symbols. In \cref{s:conv}, we state the main technical lemma, \cref{l:laplaceconverge}, that the Laplace transforms of the open ASEP height functions $\{h_{u,v,w,r}^{(N)}(\cdot)-h_{u,v,w,r}^{(N)}(0)\}_{N\in\mathbb{N}}$ converge pointwise to the Laplace transform of the open KPZ stationary measure when $u+v>0$ and $w,r>0$. The remaining parts of \cref{s:conv} demonstrate that the scaling limit of the Askey-Wilson measures produces the continuous dual Hahn process under the scaling assumptions in \cref{d:scalenew}. In \cref{s:p1} we prove \cref{l:laplaceconverge}. In \cref{s:proof}, we combine these results to prove \cref{t:1}.

\subsection{Acknowledgements} The author thanks Ivan Corwin for suggesting the problem, and Ivan Corwin and Alisa Knizel for helpful discussions about their work. The author also thanks W{\l}odek Bryc and Yizao Wang for helpful comments on a draft. The author's research was supported by Ivan Corwin’s NSF grant DMS-1811143 and the Fernholz Foundation’s “Summer Minerva Fellows” program.

\section{Background} \label{s:def}
In \cref{ss:AW} and \cref{ss:CDH} we will characterize the Askey-Wilson process and the continuous dual Hahn process, which are used to describe the stationary measures of the open ASEP and open KPZ equation, respectively. In \cref{ss:asymptote} we will quote asymptotic results on $q$-Pochhammer symbols, proved in \cite{CK21}, which will be necessary for our argument. In \cref{s:laplace}, we formally state the result that obtains weak convergence of random variables from pointwise convergence of their Laplace transforms.
\subsection{General Notation} Throughout the article, we will use the notation $(a;q)_{k}:=\prod_{j=0}^{k-1}(1-aq^{j})$ for $q$-Pochhammer symbols, and the shorthand $(a_{1},...,a_{n};q)_{k}:=(a_{1},q)_{k}\cdot\cdot\cdot(a_{n};q)_{k}$. We will also often use $(a)_{k}$ and $(a_{1},...,a_{n})_{k}$ to refer to the same quantities, removing $q$ from the notation, since for the entire paper we will set $q=\exp{\left(\frac{-2}{\sqrt{N}}\right)}$. We will likewise use $[a]_{k}:=\prod_{j=0}^{k-1}(a+j)$ to denote the standard Pochhammer symbol and a similar shorthand for products, $[a_{1},...,a_{n}]_{k}:=[a_{1}]_{k}\cdot\cdot\cdot [a_{n}]_{k}$. We will use the notation $\Gamma(x_{1},...,x_{n}):=\Gamma(x_{1})\cdot\cdot\cdot \Gamma(x_{n})$ to denote a product of Gamma functions.
\subsection{The Askey-Wilson Process}\label{ss:AW}
The stationary measure for the open ASEP can be described in terms of a stochastic process called the Askey-Wilson process. This process first appeared in \cite{BWe17}, which formulated the probability generating function of the open ASEP in terms of the Askey-Wilson process $\mathbb{Y}_{t}$, see \cref{d:awprocess}, 

$$\left\langle \prod_{j=1}^{N}t_{j}^{\tau_{j}}\right\rangle_{N}=\frac{\mathbb{E}\left[\prod_{j=1}^{N}(1+t_{j}+2\sqrt{t_{j}}\mathbb{Y}_{t_{j}})\right]}{2^{N}\mathbb{E}\left[(1+\mathbb{Y}_{1})^{N}\right]}.$$

In this paper, we will use the Askey-Wilson formula to characterize the Laplace transform of the open ASEP height function. We recall the definition of the diffusively scaled height function, $h^{(N)}(x)$ from \cref{d:scalenew}. We can write its multi-point Laplace transform as

\begin{equation}\label{e:laplace}
\begin{aligned}
    \phi^{(N)}_{u,v,w,r}(\vec{c},\vec{x}) & =\left\langle e^{-\sum_{k=1}^{d}c_{k}\left(h^{(N)}_{u,v,w,r}(x_{k}^{(N)})-h^{(N)}_{u,v,w,r}(0)\right)}\right\rangle_{N} \\ &  = \frac{\mathbb{E}\left[\prod_{k=1}^{d+1}\left(\cosh(s_{k}^{(N)})+\mathbb{Y}_{e^{-2s_{k}^{(N)}}}\right)^{n_{k}-n_{k-1}}\right]}{\mathbb{E}\left[(1+\mathbb{Y}_{1})^{N}\right]}.
\end{aligned}
\end{equation}

In this expression, as in the rest of the paper, we use the notation
\begin{equation}\label{e:assumptions}
\begin{aligned}
\vec{x} & := (x_{0},...,x_{d+1}), \text{   where   } 0=x_{0}<x_{1}<...<x_{d}\leq x_{d+1} = 1, \\ \vec{c} & := (c_{1},...,c_{d}), \text{   where   } c_{1},...,c_{d}>0,  \\ s_{k} & = c_{k} + ...+ c_{d}, \text{  and   } s_{d+1}=0.
\end{aligned}
\end{equation}
We will also use the notation $n_{k} = \lfloor Nx_{k}\rfloor$, $s_{k}^{(N)}=N^{-\frac{1}{2}}s_{k}$ and $x_{k}^{(N)}=N^{-1}\lfloor N x_{k}\rfloor$. The rest of this subsection is devoted to defining the Askey-Wilson process, following the definitions in \cite{BWe17}.

\begin{defn}
\label{d:awprocess} 
The Askey-Wilson process $\{\mathbb{Y}_{t}\}_{t\geq 0}$ is the Markov process whose marginal and transitional distributions on $V\in\mathcal{B}(\mathbb{R})$ are given by 
\begin{equation}\label{e:singlepoint}
\begin{aligned}
\pi_{t}(V) & :=AW\left(V;A\sqrt{t},B\sqrt{t},\frac{C}{\sqrt{t}},\frac{D}{\sqrt{t}}\right), \\
     \pi_{s,t}(x,V) & :=AW\left(V;A\sqrt{t},B\sqrt{t},\sqrt{\frac{s}{t}}\left(x+\sqrt{x^{2}-1}\right),\sqrt{\frac{s}{t}}\left(x-\sqrt{x^{2}-1}\right)\right).
\end{aligned}
\end{equation} In these expressions, the functions $AW(V;a,b,c,d)$ are the Askey-Wilson measures. These measures have a continuous support and a discrete support, which we denote $S^{c}(a,b,c,d,q)$ and $S^{d}(a,b,c,d,q)$, respectively. This measure is given by
\begin{align*}
    AW(V;a,b,c,d,q) & := \int_{V}AW^{c}(x;a,b,c,d,q)dx + AW^{d}(V;a,b,c,d,q), \\ AW^{c}(x;a,b,c,d,q) & :=\mathbf{1}_{|x|<1}\frac{(q,ab,ac,ad,bc,bd,cd)_{\infty}}{2\pi (abcd)_{\infty}\sqrt{1-x^{2}}}\bigg|\frac{(e^{2i\theta_{x}})_{\infty}}{(ae^{i\theta_{x}},be^{i\theta_{x}}, ce^{i\theta_{x}},de^{i\theta_{x}})_{\infty}}\bigg|^{2}, \\ AW^{d}(V;a,b,c,d,q) & :=\sum_{y\in V\cap S^{d}(a,b,c,d,q)}AW^{d}(y;a,b,c,d,q).
\end{align*} 
We refer to $AW^{c}(x;a,b,c,d,q)$ as the continuous density at $x$ and $AW^{d}(y;a,b,c,d,q)$ as the discrete mass at $y$. For an open ASEP of lattice size $N$, with parameters scaled according to \cref{d:scalenew}, we will denote the corresponding marginal and transition laws by $\pi_{s}^{(N)}(\cdot)$ and $\pi_{s,t}^{(N)}(\cdot,\cdot)$. 
\end{defn}
The discrete support of the Askey-Wilson measure  depends on the parameters $\chi\in\{a,b,c,d\}$. Atoms exist whenever one of these parameters satisfies $|\chi|>1$. If one or more of the parameters satisfies this condition, then each $\chi$ with $|\chi|>1$ will have an associated collection of atoms, which we will denote by $\{y_{j}^{\chi}\}_{j\in J}$. The index set $J$ consists of all $j$ such that $|q^{j}\chi|>1$. The atoms $y_{j}^{\chi}$ will have locations given by $y_{j}^{\chi}:=\frac{1}{2}\left(q^{j}\chi +\frac{1}{q^{j}\chi}\right).$ In the case where $\chi=a$, the discrete mass at each of the $y_{j}^{a}$ is given by

\begin{equation}\label{e:atomwght}
\begin{aligned}
AW^{d}(y_{0}^{a};a,b,c,d,q) & :=\frac{(a^{-2}, bc,bd,cd)_{\infty}}{(b/a, c/a, d/a, abcd)_{\infty}},\\
AW^{d}(y_{j}^{a};a,b,c,d,q) & :=AW^{d}(y_{0}^{a};a,b,c,d,q)\frac{(a^{2},ab,ac,ad)_{j}(1-a^{2}q^{2j})q^{j}}{(q,qa/b,qa/c,qa/d)_{j}(1-a^{2})(abcd)^{j}}. 
\end{aligned}
\end{equation}
When $\chi\neq a$, the expressions for the weights are the same, with $\chi$ and $a$ swapped throughout the expression. 
We use the notation $y_{j}^{\chi}(t)$ to represent an atom of the measure $\pi_{t}(V)$.  
\subsubsection{Atoms in the Marginal Law}\label{s:singleptatoms}
In \cref{d:awprocess}, the marginal distribution $\pi_{t}(V)$ takes the form of an Askey-Wilson measure $AW(V;a,b,c,d,q)$ with parameters 
\begin{align*}a = A\sqrt{t}=q^{v}\sqrt{t}, & & b = B\sqrt{t}=-q^{w}\sqrt{t}, & & c = C/\sqrt{t}=q^{u}/\sqrt{t},& & d=  D/\sqrt{t}=-q^{r}/\sqrt{t}.\end{align*} In order to simplify the asymptotic computations, we will deal with a rescaled process, which is indexed by time variable $q^{t}$, so that \begin{align*}a=q^{v+\frac{t}{2}},& & b=-q^{w+\frac{t}{2}}, & & c=q^{u-\frac{t}{2}},& & d=-q^{r-\frac{t}{2}}.\end{align*}
Atoms of the Askey-Wilson measure arise whenever $\chi\in\{a,b,c,d\}$ has the property that $|\chi|>1$. This can occur in one of four cases:
\begin{align*}
    \text{(1) }2v+t<0, & & \text{(2) }2w+t<0, & & \text{(3) }2u-t<0, & & \text{(4) }2r-t<0. 
\end{align*}
We recover the case of \cite{CK21} by setting $w=r=1$. In this setting, the domain in $t$ on which convergence occurs is restricted to an interval $(-2,2)$ in order to prevent atoms associated to the parameters $b$ in case $(2)$ and to $d$ in case $(4)$. When $w$ and $r$ are allowed to vary, there is a comparable domain in which there are no atoms arising from cases $(2)$ and $(4)$. For any $t\in(-2w,2r)$ we will avoid all such atoms. In particular, since $w,r>0$, $(-2w,2r)$ is an open interval. As a reminder for why this is necessary, note that establishing convergence in distribution through the Laplace transform, as described in \cref{t:laplace}, requires convergence on an open set. 

Likewise, atoms exist in cases $(1)$ and $(3)$ when $v$ and $u$ are respectively less than $1$. Since we constrain $AC=q^{u+v}<1$, we know that $u+v>0$ and therefore only one of the conditions $(1)$ or $(3)$ can hold at any given time. This fact will be important to calculations in \cref{s:conv}.

\subsubsection{Atoms and Time Reversal}\label{s:atoms}
In the discussion of the asymptotics of the transition probabilities, we center and rescale the Askey-Wilson Process. If the original Askey-Wilson process is $\mathbb{Y}$, characterized by $q,A,B,C,D,$ then we define the tangent process $$\widehat{\mathbb{Y}}^{(N)}_{s}=2N(1-\mathbb{Y}_{q^{s}}).$$ This parameterization involves a time reversal of $\mathbb{Y}$, and the transition probabilities of this process are related to the transition probabilities of the original process by a conjugation, which we will describe below. Though we will not take this approach, it is relevant to note that it is possible to avoid performing a time reversal by taking the time transformation $t\mapsto e^{2t/\sqrt{N}}$ instead of $t\mapsto q^{t}$. Other papers which employ the Askey-Wilson process, for instance~\cite{WY24}, take this approach.

We denote the measure associated to $\widehat{\mathbb{Y}}_{s}$ by $\widehat{\pi}_{s}^{(N)}$. As before, this measure has a discrete support and a continuous support, given by 
\begin{align*}\widehat{S}_{s}^{(N),c}=\widehat{S}^{(N),c}:=[0,4N], & & \widehat{S}_{s}^{(N),d}:=\left\{y\in\mathbb{R} : 1-\frac{y}{2N}\in S^{d}_{q^{s}}\right\},
\end{align*} where $S^{d}_{q^{s}}$ is the collection of atoms in $S^{d}(q^{v+s/2},-q^{1+s/2},q^{u-s/2},q^{1-s/2},q)$, this set is described in \cref{ss:AW}. In the continuous component of the support, we will additionally emphasize the dependence on the parameters $w,v,u,r$ by writing $\widehat{\pi}_{s;w,v,u,r}^{(N)}$. In this notation, we can define $\widehat{\pi}_{t}$ and $\widehat{\pi}_{s,t}$ in terms of $\pi_{t}$ and $\pi_{s,t}$ via the following transformations,
\begin{align*}
    \widehat{\pi}_{t;v,w,u,r}^{(N),c}(y)& :=\frac{1}{2N}\pi_{q^{t}}^{(N),c}\left(1-\frac{y}{2N}\right), & & y\in \widehat{S}_{t}^{(N),c},\\ \widehat{\pi}_{t}^{(N),d}(y) & := \pi_{q^{t}}^{(N),d}\left(1-\frac{y}{2N}\right), & & y\in \widehat{S}^{(N),d}_{t}.
\end{align*}
\begin{align*}
    \widehat{\pi}_{s,t;v,w,u,r}^{(N),c,c}(x,y) & :=\frac{1}{2N}\pi_{q^{t},q^{s}}^{(N),c,c}\left(1-\frac{y}{2N},1-\frac{x}{2N}\right)\frac{\pi_{q^{t}}^{(N),c}(1-\frac{y}{2N})}{\pi_{q^{s}}^{(N),c}(1-\frac{x}{2N})},  & & x\in\widehat{S}_{s}^{(N),c}, y\in\widehat{S}_{t}^{(N),c},\\
    \widehat{\pi}_{s,t;v,w,u,r}^{(N),d,c}(x,y) & :=\frac{1}{2N}\pi_{q^{t},q^{s}}^{(N),c,d}\left(1-\frac{y}{2N},1-\frac{x}{2N}\right)\frac{\pi_{q^{t}}^{(N),c}(1-\frac{y}{2N})}{\pi_{q^{s}}^{(N),d}(1-\frac{x}{2N})}, & & x\in\widehat{S}_{s}^{(N),d}, y\in\widehat{S}_{t}^{(N),c}, \\
    \widehat{\pi}_{s,t;v,w,u,r}^{(N),c,d}(x,y) & :=\pi_{q^{t},q^{s}}^{(N),d,c}\left(1-\frac{y}{2N},1-\frac{x}{2N}\right)\frac{\pi_{q^{t}}^{(N),d}(1-\frac{y}{2N})}{\pi_{q^{s}}^{(N),c}(1-\frac{x}{2N})}, 
    & & x\in\widehat{S}_{s}^{(N),c}, y\in\widehat{S}_{t}^{(N),d}, \\ 
    \widehat{\pi}_{s,t}^{(N),d,d}(x,y) & :=\pi_{q^{t},q^{s}}^{(N),d,d}\left(1-\frac{y}{2N},1-\frac{x}{2N}\right)\frac{\pi_{q^{t}}^{(N),d}(1-\frac{y}{2N})}{\pi_{q^{s}}^{(N),d}(1-\frac{x}{2N})}, & & x\in\widehat{S}_{s}^{(N),d}, y\in\widehat{S}_{t}^{(N),d}.
\end{align*}

The discrete support of $\widehat{\pi}_{q^{t},q^{s}}^{(N),d,d}\left(1-\frac{y}{2N},\cdot\right)$ is not equal to the set $\widehat{S}_{s}^{(N),d}$, and likewise for $\widehat{\pi}_{q^{t},q^{s}}\left(\cdot, 1-\frac{x}{2N}\right)$. These transition probabilities are set to be zero when $x$ and $y$ are not located, respectively, in $\widehat{S}_{s}^{(N),d}$ or $\widehat{S}_{t}^{(N),d}$. 

On its continuous support, $\widehat{\pi}^{(N),c}_{t;v,w,u,r}(\cdot)$ has the following symmetry,
\begin{align}\label{e:sym}
    \widehat{\pi}_{t;v,w,u,r}^{(N),c}(y) = \widehat{\pi}_{t;w,v,r,u}^{(N),c}(4N-y).
\end{align}

The Askey-Wilson has a limiting process called the continuous dual Hahn process, which was defined in \cite{CK21} in order to characterize the stationary measure of the open KPZ equation.

\subsection{The Continuous Dual Hahn Process}\label{ss:CDH} 
In this section, we give the definition of the continuous dual Hahn process. The definition in this paper is restricted to the time interval $t\in [0,C_{u,r})$ where $C_{u}=2r$ if $u\leq 0$ or $u\geq 1$ and $C_{u} = \min{\{2u,2r\}}$ if $u\in(0,1)$. It is possible to define the continuous dual Hahn process for all $t\in\mathbb{R}$, but the more general formulas are significantly more complicated, and not necessary for the arguments in this paper. The definition of the continuous dual Hahn process for $u,v>0$ also appears in \cite{Br12}. In definitions which follow, $\mathfrak{p}_{t}(\cdot)$ will not be the density function of a probability distribution, but instead of a positive distribution of infinite mass, while $\mathfrak{p}_{s,t}(x,\cdot)$ will correspond to a probability distribution.  We follow the definition from \cite[Section 6.2]{CK21}.
\begin{defn}\label{d:cdhproc} The \textit{continuous dual Hahn process} is the process with marginal and transition measures given by 
\begin{align}
    \mathfrak{p}_{t}(V) :=\int_{V\cap S_{t}^{c}}\mathfrak{p}_{t}^{c}(x)dx +\sum_{x\in S_{t}^{d}\cap V}\mathfrak{p}_{t}^{d}(x),  & & \mathfrak{p}_{s,t}(x,V) := \textsc{CDH}(V; a(x), b(x), c(x); \mathcal{F}),
\end{align} where $\mathfrak{p}^{c}_{t}(\cdot)$ is the marginal continuous dual Hahn law on the continuous support of the process, and $\mathfrak{p}^{d}_{t}(\cdot)$ is likewise the discrete mass defined on the discrete support of the process, and $\textsc{CDH}(V;a(x),b(x),c(x);\mathcal{F})$ is the continuous dual Hahn measure (See \cref{d:cdh}) for $x\in\mathcal{F}$ and $a(x),b(x),c(x)$ constants depending on $x$. These measures are supported on a set with a continuous component, $S^{c}_{t}$, and a discrete component, $S^{d}_{t}$, which is further broken up into two parts. We define these supports as
\begin{align*}
    S^{c}_{t}:=\mathbb{R}_{>0}, & & S^{d,u}_{t}:=\bigcup_{j=0}^{\lfloor -u+\frac{t}{2}\rfloor} \{x_{j}^{u}(t)\}, & & S_{t}^{d,v}:=\bigcup_{j=0}^{\lfloor -v-\frac{t}{2} \rfloor}\{x_{j}^{v}(t)\}. 
\end{align*} The first of the discrete collections is nonempty only when $u-t/2<0$, and the latter is nonempty only when $v+t/2<0$. The atoms in these sets are defined as
\begin{align*}
    x_{j}^{u}(t):= -4\left(u+j-\frac{t}{2}\right)^{2}, & & x_{j}^{v}(t):= -4\left(v+j+\frac{t}{2}\right)^{2}.
\end{align*}
On the continuous support, we define
\begin{align}\label{e:cdhmc}
\mathfrak{p}_{t}^{c}(x):=\frac{(u+v+1)(u+v)}{8\pi}\cdot\frac{|\Gamma\left(v+\frac{t}{2}+i\frac{\sqrt{x}}{2}\right)\Gamma\left(u-\frac{t}{2}+i\frac{\sqrt{x}}{2}\right)|^{2}}{\sqrt{x}|\Gamma(i\sqrt{x})|^{2}}\mathbf{1}_{x>0}.
\end{align} On the discrete support, we define
\begin{align}\label{e:cdhUmd}
    \mathfrak{p}_{t}^{d}(x_{j}^{u}(t)):=\frac{\Gamma(v-u+t)\Gamma(v+u+2)}{\Gamma(-2u+t)}\cdot\frac{\left(u+j-\frac{t}{2}\right)[2u-t,v+u]_{j}}{\left(u-\frac{t}{2}\right)[1,1-v+u-t]_{j}},
\end{align}
\begin{align}\label{e:cdhVmd}
    \mathfrak{p}_{t}^{d}(x_{j}^{v}(t)):=\frac{\Gamma(u-v-t)\Gamma(v+u+2)}{\Gamma(-2v-t)}\cdot\frac{\left(v+j+\frac{t}{2}\right)[2v+t,v+u]_{j}}{\left(v+\frac{t}{2}\right)[1,1-v+u-t]_{j}}.
\end{align}
To define the transition densities, we use $\mathcal{S}_{t}:=S_{t}^{c}\cup S_{t}^{d,u}\cup S_{t}^{d,v}$ to denote the union of the continuous and discrete supports at time $t$. 
In the definition which follows, we use the notation $\textsc{CDH}(V;a,b,c;\mathcal{F})$ to denote the continuous dual Hahn measure, which is introduced in \cref{d:cdh}. For $x\in S_{s}^{c}$, we define 
\begin{align}
    \mathfrak{p}_{s,t}(x,V):=\textsc{CDH}\left(V;u-\frac{t}{2}, \frac{t-s}{2}+i\frac{\sqrt{x}}{2}, \frac{t-s}{2} - i\frac{\sqrt{x}}{2}; \mathcal{S}_{t}\right).
\end{align} When $V$ is contained in the discrete support, we use the notation $\mathfrak{p}_{s,t}^{c,d}$, and likewise when $V$ is contained in the continuous support we write $\mathfrak{p}_{s,t}^{c,d}$. If $u-\frac{s}{2}<0$ (and consequently, $S_{s}^{d}=S_{s}^{d,u}$), then for $j\in [[0, \lfloor -u +\frac{s}{2}\rfloor]]$,
    \begin{align}
        \mathfrak{p}_{s,t}(x_{j}^{u}(s),V):=\textsc{CDH}\left(V;u-\frac{t}{2}, -u+\frac{t}{2}-j, u+\frac{t}{2} -s+j;\mathcal{S}_{t}\right).
    \end{align} If $v+\frac{s}{2}<0$ (and consequently, $S^{d}_{s}=S_{s}^{d,v}$, then for $j\in[[0,\lfloor -v-\frac{s}{2}\rfloor ]]$, 
    \begin{align}
        \mathfrak{p}_{s,t}(x_{j}^{v}(s),V):=\textsc{CDH}\left(V;v+j+\frac{t}{2}, \frac{t}{2}-s-v-j,u-\frac{t}{2};\mathcal{S}_{t}\right).
    \end{align}
    In both of the preceding cases, we use the notation $\mathfrak{p}_{s,t}^{d,d}$ and $\mathfrak{p}_{s,t}^{d,c}$ to denote the transition to the discrete and continuous supports, respectively. In all other cases, we define $\mathfrak{p}_{s,t}$ to be zero.
\end{defn}
The next definition, of $\textsc{CDH}^{c}(x;a,b,c)$ and $\textsc{CDH}^{d}(x_{j};a,b,c;\mathcal{F})$ follows \cite[Definitions 6.3 and 6.4]{CK21}.
\begin{defn}\label{d:cdh} Assume $a,b,c\in\mathbb{R}$ with $a+b,a+c>0$ or $a\in\mathbb{R}$ and $b=\overline{c}$ with $Re(b)=Re(c)>0$. For such $a,b,c$ we define  
\begin{align*}
    \textsc{CDH}^{c}(x;a,b,c):= \frac{1}{8\pi}\cdot\frac{\bigg|\Gamma\left(a+i\frac{\sqrt{x}}{2}\right)\Gamma\left(b+i\frac{\sqrt{x}}{2}\right)\Gamma\left(c+i\frac{\sqrt{x}}{2}\right)\bigg|^{2}}{\Gamma(a+b)\Gamma(a+c)\Gamma(b+c)\sqrt{x}|\Gamma(i\sqrt{x})|^{2}}\mathbf{1}_{x>0}.
\end{align*} Similarly, in the setting where $a<0$ and $\mathcal{F}$ is a finite subset of $\mathbb{R}$ with size $|\mathcal{F}|=\lfloor-a\rfloor +1$ and elements $x_{0}<x_{1}<\cdot\cdot\cdot<x_{\lfloor -a\rfloor}$, then assume that $b,c,$ and $j$ satisfy one of the following conditions:
\begin{enumerate}
    \item $b=\overline{c}$ with $Re(b),Re(c)>0$ and $j\in [[0,\lfloor -a\rfloor]]$
    \item $b,c\in\mathbb{R}$ with $a+b,a+c>0$ and $j\in [[0,\lfloor -a\rfloor]]$
    \item $b,c\in\mathbb{R}$ and $b,b+c,c-a>0$ with $a+b=-k$ for $k\in\mathbb{Z}_{\geq 0}$ and $j\in[[0,k]]$. 
\end{enumerate} When one of these conditions is satisfied, we define the discrete mass at $x_{j}$ to be
\begin{align*}
    \textsc{CDH}^{d}(x_{j};a,b,c;\mathcal{F}):=\frac{[2a,a+b,a+c]_{j}(a+j)\Gamma(b-a)\Gamma(c-a)}{[1,a-b+1,a-c+1]_{j}a\Gamma(-2a)\Gamma(b+c)}(-1)^{j}.
\end{align*} On all other arguments, $\textsc{CDH}^{d}(\cdot;a,b,c;\mathcal{F})$ is zero. Finally, we define the continuous dual Hahn measure $\textsc{CDH}(V;a,b,c;\mathcal{F})$ in the following three cases:
\begin{itemize}
    \item For $a\geq 0$; $\mathcal{F}=\emptyset$ and either $b=\overline{c}$ with $Re(b),Re(c)>0$ or $b,c\in\mathbb{R}_{>0}$, $$\textsc{CDH}(V;a,b,c;\mathcal{F}):=\int_{V}\textsc{CDH}^{c}(x;a,b,c)dx.$$
    \item For $a<0$; $\mathcal{F}$ a finite subset of $\mathbb{R}$ with size $|\mathcal{F}|=\lfloor -a\rfloor + 1$ and elements $x_{0}<x_{1}<\cdot\cdot\cdot <x_{\lfloor -a\rfloor}$, and either $b=\overline{c}$ with $Re(b),Re(c)>0$ or $b,c\in\mathbb{R}$ with $a+b,a+c>0$, $$\textsc{CDH}(V;a,b,c;\mathcal{F}):=\int_{V}\textsc{CDH}^{c}(x;a,b,c)dx + \sum_{x\in\mathcal{F}\cap V}\textsc{CDH}^{d}(x;a,b,c;\mathcal{F}).$$
    \item For $a<0$; $\mathcal{F}$ a finite subset of $\mathbb{R}$ with size $|\mathcal{F}|=\lfloor -a\rfloor + 1$ and elements $x_{0}<x_{1}<\cdot\cdot\cdot<x_{\lfloor -a\rfloor}$ and $b,c\in\mathbb{R}$ with $b,b+c,c-a>0$ and $a+b=-k$ for $k\in\mathbb{Z}_{\geq 0}$ $$\textsc{CDH(}V;a,b,c;\mathcal{F}):=\sum_{x\in\mathcal{F}\cap V}\textsc{CDH}^{d}(x;a,b,c;\mathcal{F}).$$
\end{itemize}

\end{defn}

\subsection{Notation for Asymptotics}\label{ss:asymptote}
We will use the following notation to simplify the expressions for the logarithms of the $q$-Pochhammer symbols.
\begin{align}\label{e:aplusminus}
    \mathcal{A}^{+}(\kappa,z) = -\frac{\pi^{2}}{6\kappa}-\left(z-\frac{1}{2}\right)\log{(\kappa)} - \log{\left(\frac{\Gamma(z)}{\sqrt{2\pi}}\right)}, & &  \mathcal{A}^{-}(\kappa,z)= \frac{\pi^{2}}{12\kappa} - \left(z-\frac{1}{2}\right)\log{(2)}.
\end{align}
The following lemma will also be crucial to our asymptotic analysis.
\begin{lemma}[Proposition 2.3~\cite{CK21}]\label{p:asymp}

For $\kappa\in(0,1)$, $z\in\mathbb{C}$, and $q=\exp{(-\kappa)}$, and $m\in\mathbb{Z}_{>0}$,
$$\log{\left((q^{z};q)_{\infty}\right)}=\mathcal{A}^{+}(\kappa,z)-\sum_{n=1}^{m-1}\frac{B_{n+1}(z)B_{n}}{n(n+1)!}\kappa^{n}+e^{+}_{m}(\kappa,z),$$
    $$\log{\left((-q^{z};q)_{\infty}\right)}=\mathcal{A}^{-}(\kappa,z)-\sum_{n=1}^{m-1}(2^{n}-1)\frac{B_{n+1}(z)B_{n}}{n(n+1)!}\kappa^{n}+e^{-}_{m}(\kappa,z),$$
where $B_{n}(z)$ and $B_{n}$ denote Bernoulli polynomials and Bernoulli numbers, respectively. For any $\alpha\in (0,\pi)$, $\epsilon\in(0,1/2)$ and any $b\in (m-1,m)$ there exist $C,\kappa_{0}>0$ such that for all $\kappa\in(0,\kappa_{0})$ and all $z\in\mathbb{C}$ with $|Im(z)|<\frac{\alpha}{\kappa}$, 
\begin{align}\label{e:errorbd}|e_{1}^{\pm}(\kappa,z)|\leq C(\kappa(1+|z|)^{2}+\kappa^{b}(1+|z|)^{1+2b+\epsilon}).\end{align} The bound on $e^{+}_{m}(\kappa,z)$ further holds with the condition $|Im(z)|<\frac{\alpha}{\kappa}$ replaced by $|Im(z)|<\frac{2\alpha}{\kappa}$. Furthermore, for any $r>0$, $\epsilon\in (0,1/2)$ and $b\in (m-1,m)$ there exist $C,\kappa_{0}>0$ such that for all $\kappa\in(0,\kappa_{0})$ and all $z\in\mathbb{C}$ with dist$(\text{Re}(z),\mathbb{Z}_{\leq 0})>r$, \eqref{e:errorbd} continues to hold. 
\end{lemma} 
We also use the following notation as shorthand \begin{align}\label{e:aplusminusN}\mathcal{A}^{\pm}_{N}\left(z_{1},...,z_{k}\right)=\sum_{i=1}^{k}\mathcal{A}^{\pm}\left(\kappa,z_{i}\right), & & E^{\pm}_{N}\left(z_{1},...,z_{k}\right)=\sum_{i=1}^{k}e_{1}^{\pm}\left(\kappa,z_{i}\right).\end{align}

\subsection{Convergence of the Laplace transform}\label{s:laplace}

Finally, we state the result about pointwise convergence of Laplace transforms, which we will use to obtain our weak convergence.

\begin{prop}[Section 5.14, page 378, (5.14.8) ~\cite{hoff}]\label{t:laplace}  
    For a collection of random vectors of length $d$, $\{\vec{X}^{(i)}\}_{i\in\mathbb{N}}$, with $\vec{X}^{(i)}:=(X^{(i)}_1,\ldots, X^{(i)}_d)$ and $\vec{z}\in\mathbb{R}^{d}$ defined so that $\vec{z}:=(z_1,\ldots, z_d)$, we define the collection of functions $\{L_{i}(\vec{z})\}_{i\in\mathbb{N}}$ to be the Laplace transforms of each of the $\vec{X}^{(i)}$, respectively, $$L_{i}(\vec{z})=\mathbb{E}\left[\exp{\left(\sum_{j=1}^{d}z_{j}X_{j}^{(i)}\right)}\right].$$ Suppose that $L_{i}(\vec{z})$ are all finite and converge pointwise to a function $L(\vec{z})$ for all $\vec{z}$ in an open set of $\mathbb{R}^{d}$. If $L(\vec{z})$ is the Laplace transform of a random variable, $\vec{Y}=(Y_{1},...,Y_{d})$, then the sequence $\vec{X}^{(i)}$ converges in distribution to $\vec{Y}.$
\end{prop}
\section{Convergence of the Askey Wilson Tangent Laws to the Continuous Dual Hahn Laws}\label{s:conv}
The purpose of this section and the next is to demonstrate that the sequence of Laplace transforms of the Askey-Wilson process, under \cref{d:scalenew}, converges pointwise to the Laplace transform of the open KPZ stationary measure.

\begin{prop}\label{l:laplaceconverge}
    The multi-point Laplace transforms $\varphi^{(N)}_{u,v,w,r}(\vec{c},\vec{x})$, with $\vec{c},\vec{x}$ as in \eqref{e:assumptions}, of the sequence of open ASEP height functions $\{h_{u,v,w,r}^{(N)}(\cdot)\}_{N\in\mathbb{N}}$, scaled as in \cref{d:scalenew}, are represented in terms of the Askey-Wilson process $\mathbb{Y}_{q^{s}}$ as
\begin{equation}
\begin{aligned}
    \varphi^{(N)}_{u,v,w,r}(\vec{c},\vec{x}) &  = \frac{\mathbb{E}\left[\prod_{k=1}^{d+1}\left(\cosh(s_{k}^{(N)})+\mathbb{Y}_{e^{-2s_{k}^{(N)}}}\right)^{n_{k}-n_{k-1}}\right]}{\mathbb{E}\left[(1+\mathbb{Y}_{1})^{N}\right]}.
\end{aligned}
\end{equation} For all $\vec{c}\in (0,C_{u,r})^{d}$, defined so that 
\begin{align*}C_{u,r} & :=\begin{cases} 2 & u \leq 0, r\leq 0, \\ 2r & r>0, u\leq 0, \\ \min{\{2r,2u\}} & u,r> 0 ,
\end{cases}
\end{align*} for $u+v>0$ and $w,r>0$, we have pointwise convergence $\lim_{N\to\infty}\phi^{(N)}_{u,v,w,r}(\vec c, \vec x) = \varphi_{u,v}(\vec c, \vec x)$. The function $\varphi_{u,v}(\vec c,\vec x)$ is defined
\begin{align}\label{e:laplacecharacterization}
    \varphi_{u,v}(\vec{c},\vec{x}) & := \frac{\mathbb{E}\left[\exp{\left(\frac{1}{4}\sum_{k=1}^{d+1}(s_{k}^{2}-\mathbb{T}_{s_{k}})(x_{k}-x_{k-1})\right)}\right]}{\mathbb{E}\left[\exp{\left(-\frac{1}{4}\mathbb{T}_{0}\right)}\right]},
\end{align} where $\mathbb{T}_{t}$ is the continuous dual Hahn process.
Moreover, $\varphi_{u,v}(\vec{c},\vec{x})$ is the Laplace transform of the stationary measure of the open KPZ equation constructed in \cite{CK21}. 
\end{prop}

The first step in the proof of \cref{l:laplaceconverge} is to demonstrate the convergence of the Askey-Wilson process measures to those of the continuous dual Hahn process, which we do in the rest of this section. In \cref{s:p1} we will complete the proof of \cref{l:laplaceconverge}. Finally, in \cref{s:proof}, we give the proof of \cref{t:1}. 

We will use the notation $x_{j}^{\theta}(t)$ and $y_{j}^{(N),\theta}$, respectively, for the locations of the atoms of the rescaled continuous dual Hahn process and Askey-Wilson process associated to $\theta \in \{v,u,w,r\}$ (with correspond to the atoms arising from, respectively, $\{q^{v+t/2},q^{u-t/2},-q^{w+t/2},-q^{r-t/2}\}$).

\begin{lemma}\label{l:singlepoint} Under the scaling of the Askey-Wilson measure when $q=\exp{\left(\frac{-2}{\sqrt{N}}\right)}$, the following convergence results hold.
\begin{enumerate}
    \item For all $y\in \widehat{S}^{c}$ and $t\in\mathbb{R}_{>0}$, for all $\eta>1$ there exist $N_{0}\in\mathbb{N}$ and  $C,\chi>0$ so that for all $N > N_{0}$ and $y\in [0,4N]$, 
\begin{align*} \label{e:conv}
     N^{u+v}\sqrt{1-\frac{y}{4N}}\widehat{\pi}_{t}^{(N),c}(y)   = C_{u,v,w,r}\cdot \mathfrak{p}_{t}^{c}(y)\cdot e^{\textsc{Error}_{t}^{(N),c}(y)},
 \end{align*} with $$|\textsc{Error}_{t}^{(N),c}(y)|  \leq CN^{-\chi}(1+\sqrt{r})^{\eta}.$$  In particular, we have 
    \begin{align*}
        \lim_{N\to\infty}N^{u+v}\widehat{\pi}_{t}^{(N),c}(y) & = C_{u,v,w,r}\cdot \mathfrak{p}_{t}^{c}(y),
    \end{align*} where $$C_{u,v,w,r}:=\frac{\Gamma(u+v+w+r)}{\Gamma(u+v+2)\Gamma(w+r)}.$$
    \item Similarly, for $\widehat{y}_{j}^{(N),\theta}(t)\in\widehat{S}^{d}_{t}$ and $t\in(-2w,2r)$, there exist $C_{1},C_{2},\eta_{1},\eta_{2}\in\mathbb{R}_{>0}$ such that for all $N\in\mathbb{Z}_{\geq 1}$, the atoms $y_{j}^{(N),\theta}(t)$ for $\theta\in\{v,u\}$ satisfy $$N^{u+v}\widehat{\pi}_{t}^{(N),d}(\widehat{y}^{(N),v}_{j}(t))=C_{u,v,w,r}\cdot\mathfrak{p}_{t}^{d}(x_{j}^{v}(t))\cdot e^{\textsc{Error}_{t}^{(N),d}(x_{j}^{v}(t))},$$
$$N^{u+v}\widehat{\pi}_{t}^{(N),d}(\widehat{y}^{(N),u}_{j}(t))=C_{u,v,w,r}\cdot\mathfrak{p}_{t}^{d}(x_{j}^{u}(t))\cdot e^{\textsc{Error}_{t}^{(N),d}(x_{j}^{u}(t))},$$ with the exponential error terms bounded by 
$$|\textsc{Error}_{t}^{(N),d}(x_{j}^{v}(t))|\leq C_{1}N^{-\eta_{1}},$$
$$|\textsc{Error}_{t}^{(N),d}(x_{j}^{u}(t))|\leq C_{2} N^{-\eta_{2}},$$ and $C_{u,v,w,r}$ as defined above. In particular, we have 
    \begin{align*}
        \lim_{N\to\infty}N^{u+v}\widehat{\pi}_{t}^{(N),d}(\widehat{y}_{j}^{(N),\chi}(t)) = C_{u,v,w,r}\cdot \mathfrak{p}_{t}^{d}(x_{j}^{\chi}(t)).
    \end{align*}
\end{enumerate}
The limiting measure $\mathfrak{p}_{t}^{c}(y)$ and $\mathfrak{p}_{t}^{d}(y)$ are the marginal measures of the continuous dual Hahn process, see \cref{ss:CDH}.
\end{lemma}
Similarly, we can compute the limit of $\widehat{\pi}_{s,t}^{(N),c,c}(\cdot,\cdot)$ and $\widehat{\pi}_{s,t}^{(N),d,c}(\cdot,\cdot)$ under scaling. 
\begin{lemma}\label{l:transprob}With scaling conditions as in \cref{d:scalenew}, and for all $s<t$, the following results hold.
    \begin{enumerate}
        \item  For all $b,r\in\widehat{S}^{(N),c}=[0,4N]$ and $\eta>0$, there exists $N_{0}\in\mathbb{Z}_{\geq 1}$ and $C,\chi\in\mathbb{R}_{>0}$ such that for all $N>N_{0}$ 
        $$\sqrt{1-\frac{r}{4N}}\widehat{\pi}_{s,t}^{(N),c,c}(b,r)=\mathfrak{p}_{s,t}^{c,c}(b,r)\cdot e^{\textsc{Error}_{s,t}^{(N),c,c}(b,r)},$$ with an error term satisfying $$|\textsc{Error}_{s,t}^{(N),c,c}(b,r)|\leq CN^{-\chi}(1+\sqrt{b})^{\eta} + CN^{-\chi}(1+\sqrt{r})^{\eta}.$$ In particular, $$\lim_{N\to\infty}\widehat{\pi}_{s,t}^{(N),c,c}(b,r)=\mathfrak{p}_{s,t}^{c,c}(b,r).$$
        \item If $s,t\in(-2b,2d)$ and $v+\frac{s}{2}<0$, 
        so that $\widehat{S}_{s}^{(N),d}=\{\widehat{y}_{j}^{(N),v}(s)\}_{j\in [[0,\lfloor -v-\frac{s}{2}\rfloor ]]}$, 
        then for each such $j$, and all $r\in\widehat{S}^{(N),c}=[0,4N]$ there exists $N_{0}\in\mathbb{Z}_{\geq 1}$ and $C,\chi\in\mathbb{R}_{>0}$ such that for all $N>N_{0}$ $$\widehat{\pi}_{s,t}^{(N),d,c}\left(\widehat{y}_{j}^{(N),v}(s),r\right)=\mathfrak{p}_{s,t}^{d,c}(x_{j}^{v}(s),r)\cdot e^{\textsc{Error}_{s,t}^{(N),d,c}(x_{j}^{v}(s),r)},$$ with $$|\textsc{Error}_{s,t}^{(N),d,c}(x_{j}^{v}(s),r)|\leq CN^{-\chi}(1+\sqrt{r})^{\eta}.$$ In particular, $$\lim_{N\to\infty}\widehat{\pi}_{s,t}^{(N),d,c}(\widehat{y}_{j}^{(N),v}(s),r)=\mathfrak{p}_{s,t}^{d,c}(x_{k}^{v}(s),r).$$

        \item If $s,t\in(-2w,2r)$ and $u-\frac{s}{2}<0$, so that $\widehat{S}_{s}^{(N),d}=\{\widehat{y}_{j}^{(N),u}(s)\}_{j\in [[0,\lfloor -u+\frac{s}{2}\rfloor ]]}$, then for each such $j$ and for all $r\in\widehat{S}^{(N),c}=[0,4N]$,
        $$\lim_{N\to\infty}\widehat{\pi}_{s,t}^{(N),d,c}(\widehat{y}_{j}^{(N),v}(s),r)=0.$$
    \end{enumerate}
\end{lemma}
Finally, we obtain the following asymptotic result for $\widehat{\pi}_{t}^{(N),d,d}(\cdot,\cdot)$ and $\widehat{\pi}_{t}^{(N),c,d}(\cdot,\cdot)$.
\begin{lemma}\label{l:ddconv}
    With scaling conditions as defined in \cref{d:scalenew}, and for all $s,t\in[0,C_{u,r})$ with $s<t$,
    \begin{enumerate}
        \item If $v+\frac{s}{2}<0$ then the  discrete support of $\widehat{\pi}_{s,t}^{(N),d,d}(\widehat{y}_{j}^{(N),v}(s), \cdot)$ coincides with the points $\widehat{y}_{k}^{(N),v}(t)$ 
        for $k\in[[0,\lfloor -v-\frac{t}{2}\rfloor]]$ 
        and there exist $C,\chi>0$ such that for all $N>0$ and $j\in[[0,\lfloor -v-\frac{s}{2}\rfloor]],$ respectively, $k\in[[0,\lfloor -u+\frac{t}{2} \rfloor]]$
        $$\widehat{\pi}_{s,t}^{(N),d,d}(\widehat{y}_{j}^{(N),v}(s),\widehat{y}_{k}^{(N),v}(t))=\mathfrak{p}_{s,t}^{d,d}(x_{j}^{v}(s),x_{k}^{v}(t))\cdot e^{\textsc{Error}_{s,t}^{(N),d,d}(v,v)},$$  
        and $$|\textsc{Error}_{s,t}^{(N),d,d}(\cdot,\cdot)|\leq CN^{-\chi}.$$
        In particular,  $$\lim_{N\to\infty}\widehat{\pi}_{s,t}^{(N),d,d}(\widehat{y}_{j}^{(N),v}(s),\widehat{y}_{k}^{(N),v}(t))=\mathfrak{p}_{s,t}^{d,d}(x_{j}^{v}(s),x_{k}^{v}(t)).$$
        \item If $u-\frac{s}{2}<0$ then the  discrete support of $\widehat{\pi}_{s,t}^{(N),d,d}(\widehat{y}_{j}^{(N),u}(s), \cdot)$ coincides with the points $\widehat{y}_{k}^{(N),u}(t)$ 
        for $k\in[[0,j]]$
        and for all such $j,k$ and the associated atoms 
        \begin{align*}
            \widehat{\pi}_{s,t}^{(N),d,d}\left(\widehat{y}_{j}^{(N),u}(s),\widehat{y}_{k}^{(N),u}(t)\right) & = \mathfrak{p}_{s,t}^{d,d}(x_{j}^{u}(s),x_{j}^{v}(t))\cdot e^{\textsc{Error}_{s,t}^{(N),d,d}(u,u)},
        \end{align*} and 
        \begin{align*}
            |\textsc{Error}_{s,t}^{(N),d,d}(\cdot,\cdot)| & \leq CN^{-\chi}.
        \end{align*}
        In particular, $$\lim_{N\to\infty}\widehat{\pi}_{s,t}^{(N),d,d}(\widehat{y}_{j}^{(N),u}(s),\widehat{y}_{k}^{(N),u}(t))=\mathfrak{p}_{s,t}^{d,d}(x_{j}^{u}(s),x_{k}^{u}(t)).$$
        \item If $v+\frac{s}{2}<0$ and $r\in\widehat{S}^{(N),c}=[0,4N]$ then the measure $\widehat{\pi}_{s,t}^{(N),c,d}(r,\cdot)$ has no atomic part. When $u-\frac{s}{2}<0$ and $r\in\widehat{S}^{(N),c}=[0,4N]$, the measure $\widehat{\pi}_{s,t}^{(N),c,d}(r,\cdot)$ has discrete support on the atoms $\widehat{y}_{k}^{(N),u}(t)$ for $k\in[[0,\lfloor -u+\frac{t}{2}\rfloor]]$ and there exist $C,\chi\in\mathbb{R}_{>0}$ such that for all $N\in\mathbb{Z}_{\geq 1}$, $$\widehat{\pi}_{s,t}^{(N),c,d}\left(r,\widehat{y}_{k}^{(N),u}(t)\right)=\mathfrak{p}_{s,t}^{c,d}(r,x_{k}^{u}(t))\cdot e^{\textsc{Error}_{s,t}^{c,d}(r,x_{k}^{u}(t))},$$ with the error bounded $$|\textsc{Error}_{s,t}^{c,d}(r,x_{k}^{u}(t))|\leq CN^{-\chi}.$$
        In particular, 
        $$\lim_{N\to\infty}\widehat{\pi}_{s,t}^{(N),c,d}(r,\widehat{y}_{k}^{(N),u}(t))=\mathfrak{p}_{s,t}^{c,d}(r,x_{k}^{u}(t)).$$ 
    \end{enumerate}
\end{lemma}

\subsection{Convergence of Marginal Distributions}\label{s:sinconv}

This section deals with the proof of \cref{l:singlepoint}.
\begin{proof}[Proof of  \cref{l:singlepoint} (1)]
We begin by reminding ourselves of the expression for $\widehat{\pi}_{t;v,w,u,r}^{(N),c}(\cdot)$. As before, we set $q=\exp{\left(-\frac{2}{\sqrt{N}}\right)}$.  We will deal with the expression for the rescaled, time-reversed process $\widehat{\mathbb{Y}}_{t}$, and the corresponding marginal measures $\widehat{\pi}_{t;v,w,u,r}^{(N),c}(\cdot)$. 

\begin{equation}\label{e:singlecmeasure}
\begin{aligned}
    \widehat{\pi}_{t;v,w,u,r}^{(N),c}(y)&=\frac{1}{2N}\pi_{q^{t}}^{(N),c}\left(1-\frac{y}{2N}\right)  \\
    & =\frac{1}{4N\pi(q^{u+v+w+r})_{\infty}\sqrt{\frac{y}{N}}\sqrt{1 - \frac{y}{4N}}} \\ & \times \frac{(q,-q^{v+w+t}, q^{u+v}, -q^{v+r}, -q^{u+w}, q^{w+r}, -q^{u+r-t})_{\infty}|(q^{i\sqrt{N}\theta})_{\infty}|^{2}}{|(q^{v+\frac{t}{2}+i\frac{\sqrt{N}\theta}{2}},-q^{w+\frac{t}{2}+i\frac{\sqrt{N}\theta}{2}}, q^{u-\frac{t}{2}+i\frac{\sqrt{N}\theta}{2}}, -q^{r-\frac{t}{2}+i\frac{\sqrt{N}\theta}{2}})_{\infty}|^{2}} . 
\end{aligned}
\end{equation} 

We will only prove a bound in the interval $y\in[0,2N]$. The way that this bound can be applied to prove the full \cref{l:singlepoint} is by using the symmetry \eqref{e:sym} to turn this bound into the proof by first proving our desired bound for $\widehat{\pi}^{(N),c}_{t;v,w,u,r}(y)$ on the interval $y\in[0,2N]$, and then noting that the desired bound on the interval $[2N,4N]$ is significantly weaker than the bound on the interval $[0,2N]$ which arises from $\widehat{\pi}^{(N),c}_{t;w,v,r,u}(4N-y)$. Therefore, we first prove the appropriate bound on \eqref{e:singlecmeasure} in the interval $r\in[0,2N]$, and then use this identity to reach the remaining $r\in[2N,4N].$
We take the logarithm of the marginal measure,
\begin{align}\label{e:singlec}
    \log{\left(\sqrt{1-\frac{y}{4N}}\widehat{\pi}_{t;v,w,u,r}^{(N),c}(y)\right)} & =-\log{\left(4\pi N\sqrt{\frac{y}{N}}\right)} + \mathcal{A}_{N}(\sqrt{y}) + E^{\mathcal{A}}_{N}(y) + E^{\mathcal{E}}_{N}(y),
\end{align}  where, in the expression above, 
\begin{align*}
    \mathcal{A}_{N}(\sqrt{y}) & : = \mathcal{A}^{+}_{N}\left(1,u+v, w+r,i\frac{\sqrt{y}}{2}, -i\frac{\sqrt{y}}{2}\right) \\ & +\mathcal{A}^{-}_{N}\left(v+w+t,v+r,u+w,u+r - t\right) \\ & - \mathcal{A}_{N}^{+}\left(u+v+w+r,v+i\frac{\sqrt{y}}{2}+\frac{t}{2},u+i\frac{\sqrt{y}}{2}-\frac{t}{2},v-i\frac{\sqrt{y}}{2}+\frac{t}{2},u-i\frac{\sqrt{y}}{2}-\frac{t}{2}\right) \\ & - \mathcal{A}_{N}^{-}\left(w+i\frac{\sqrt{y}}{2} + \frac{t}{2},r+i\frac{\sqrt{y}}{2}-\frac{t}{2}, w-i\frac{\sqrt{y}}{2} + \frac{t}{2},r-i\frac{\sqrt{y}}{2}-\frac{t}{2}\right),
\\ E^{\mathcal{A}}_{N}(y) & :=\mathcal{A}_{N}(\sqrt{y}+E^{\theta}_{N}(y))-\mathcal{A}_{N}(\sqrt{y}), \\
E_{N}^{\mathcal{E}}(y) & :=\mathcal{E}_{N}(\sqrt{y}+E_{N}^{\theta}(y)),\\
E_{N}^{\theta}(y) & :=\sqrt{N}\arccos{\left(1-\frac{y}{2N}\right)}-\sqrt{y},\end{align*} are a concise way of containing the logarithm of the $q$-Pochhammer symbols in $\widehat{\pi}_{t;v,w,u,r}^{(N),c}(y)$. 
The notation $\mathcal{A}^{\pm}$ was defined in \eqref{e:aplusminus}, and $\mathcal{A}^{\pm}_{N}$ defined in \eqref{e:aplusminusN}, and $E_{N}^{\theta}(y)$ arises as a concise expression for the error terms in the Taylor expansion $\sqrt{N}\theta_{y}^{(N)}=\sqrt{y}+E_{N}^{\theta}(y)$.

We can simplify the first two terms on the right-hand side of \eqref{e:singlec} considerably by expanding $\mathcal{A}(\sqrt{y})$ through the definitions of $\mathcal{A}^{+}_{N}(z_{1},...,z_{k})$ and $\mathcal{A}^{-}_{N}(z_{1},...,z_{k})$ from \eqref{e:aplusminusN}, and then canceling where appropriate. The result of that calculation yields
\begin{align*}
     -\log{\left(4\pi N\sqrt{\frac{y}{N}}\right)}+\mathcal{A}_{N}(\sqrt{y}) 
     & =-\left(u+v\right)\log{(N)} \\  & +\log{\left(\frac{\Gamma\left(u+v+w+r\right)\bigg|\Gamma\left(v+\frac{i\sqrt{y}}{2}+\frac{t}{2}\right)\bigg|^{2}\bigg|\Gamma\left(u+\frac{i\sqrt{y}}{2}-\frac{t}{2}\right)\bigg|^{2}}{8\pi\sqrt{y}\Gamma(u+v)\Gamma\left(w+r\right)|\Gamma(i\sqrt{y}/2)|^{2}}\right)}.
\end{align*}
Using the definitions of $\mathfrak{p}_{t}^{c}(y)$ and $C_{u,v,w,r}$ from \cref{l:singlepoint}, we can rewrite \eqref{e:singlec} as
\begin{align*}
    \log\left(N^{u+v}\sqrt{1-\frac{y}{4N}}\widehat{\pi}_{t;v,w,u,r}^{(N),c}(y)\right) & = \log{(\mathfrak{p}_{t}^{c}(y))}+ \log{(C_{u,v,w,r})}+ E_{N}^{\mathcal{A}}(y) + E_{N}^{\mathcal{E}}(y).
\end{align*}
What remains is to check that the error terms are bounded appropriately. These arguments are identical to those in \cite{CK21}, except that the error terms may now vary in $w$ and $r$. We summarize them here for completeness. First we will quote, as a lemma, a few intermediate results from \cite{CK21}, which provide helpful bounds. 

\begin{lemma}[p. 53-54~\cite{CK21}]\label{l:ethetanbound}
    \begin{enumerate}
        \item For all $\eta>1$ there exists $N_{0}\in\mathbb{Z}_{\geq 1}$ and $C,\chi>0$ such that for all $N>N_{0}$ and $y\in\widehat{S}^{(N),c}=[0,2N]$, $$E_{N}^{\theta}(y)\log{(N)}\leq CN^{-\chi}(1+\sqrt{y})^{\eta}.$$
        \item For any fixed $z$ and for all $b\in(0,1)$ there exist $\kappa_{0}>0$ and $C>0$ such that for all $\kappa<\kappa_{0}$, $$|e_{1}^{\pm}(\kappa,z)|\leq C(\kappa+\kappa^{b})\leq C'N^{-\frac{b}{2}}.$$ 
        \item For $z_{c}(y)=c\pm i\frac{\sqrt{y}+E_{N}^{\theta}(y)}{2}$, for any $b\in(0,1)$ and  $\epsilon\in(0,1/2)$, there exists $C>0$ and $N_{0}\in\mathbb{Z}_{\geq 1}$ such that for all $N>N_{0}$ and $y\in\widehat{S}^{(N),c}=[0,2N]$, $$|e_{1}^{\pm}(\kappa,z_{c}(y))|\leq C\left(N^{-1/2}(1+\sqrt{y})^{2}+N^{-b/2}(1+\sqrt{y})^{1+2b+\epsilon}\right).$$
        \item For $z_{c}(y)=c\pm i \frac{\sqrt{y}+E_{N}^{\theta}(y)}{2}$, $z'_{c}(y)=c\pm i \frac{\sqrt{y}}{2}$ and for any $c\in\mathbb{R}$ fixed, there exists $N_{0}\in\mathbb{Z}_{\geq 1}$ and $C,\chi>0$ such that for all $N>N_{0}$ and $y\in [0,4N]$, $$|\log{\left(\Gamma(z'_{c}(y))\right)}-\log{\left(\Gamma(z_{c}(y))\right)}|\leq CN^{-\chi}(1+\sqrt{y})^{\eta}.$$
    \end{enumerate}
\end{lemma}

With these results in hand, we can state the bound on the error terms.
\begin{lemma}\label{p:ebounds}
For all $\eta>1$ there exists $N_{0}$ and $C,\chi>0$ 
 such that for all $N> N_{0}$ and $y\in[0,2N]$,  \begin{align}\label{e:ebounds}
     |E_{N}^{\mathcal{A}}(y)|,|E_{N}^{\mathcal{E}}(y)|\leq C N^{-\chi}(1+\sqrt{y})^{\eta}.
 \end{align}
\end{lemma}
\begin{proof}
There are two types of terms involved in $E^{\mathcal{E}}_{N}(y)$. We separate the terms of $E_{N}^{\mathcal{E}}(y)$ into terms which only depend on $u,v$ and $t$ (type (1)) and terms which depend on $y$ (type (2)). 

\begin{align*}
    \textbf{Type (1)} & = E_{N}^{+}\left(1, u+v, w+r\right) + E_{N}^{-}\left(v+t+w,v+r,u+w,u+r-t\right) - E_{N}^{+}\left(u+v+w+r\right), \\ \textbf{Type (2)} & = E_{N}^{+}\left(\frac{i\sqrt{N}}{2}\arccos{\left(1-\frac{y}{4N}\right)}, \frac{-i\sqrt{N}}{2}\arccos{\left(1-\frac{y}{4N}\right)}\right) \\ & - E_{N}^{+}\left(v+\frac{t}{2}+\frac{i\sqrt{N}}{2}\arccos{\left(1-\frac{y}{4N}\right)}, v+\frac{t}{2}-\frac{i\sqrt{N}}{2}\arccos{\left(1-\frac{y}{4N}\right)}\right) \\ & -E_{N}^{+}\left(u-\frac{t}{2}+\frac{i\sqrt{N}}{2}\arccos{\left(1-\frac{y}{4N}\right)},u-\frac{t}{2}-\frac{i\sqrt{N}}{2}\arccos{\left(1-\frac{y}{4N}\right)}\right) \\ & - E_{N}^{-}\left(w+\frac{t}{2} + \frac{i\sqrt{N}}{2}\arccos{\left(1-\frac{y}{4N}\right)}, w + \frac{t}{2} - \frac{i\sqrt{N}}{2}\arccos{\left(1-\frac{y}{4N}\right)}\right) \\ & - E_{N}^{-}\left(r - \frac{t}{2} + \frac{i\sqrt{N}}{2}\arccos{\left(1-\frac{y}{4N}\right)},r-\frac{t}{2}-\frac{i\sqrt{N}}{2}\arccos{\left(1-\frac{y}{4N}\right)}\right).
\end{align*}
We note that the terms of type (1) do not vary in $y$ and we can therefore apply the bound from \cref{l:ethetanbound} (2). 
The terms of type (2) come from expressions of the form $\textsc{Error}_{1}^{\pm}[\kappa,z]$ with $z$ taking the form $z_{c}(y)=c+i\left(\frac{\sqrt{y}+E_{N}^{\theta}(y)}{2}\right)$ for $c\in\mathbb{R}$. Therefore, we can directly apply \cref{l:ethetanbound} (3) to bound these terms. 

Turning our attention to $E^{\mathcal{A}}_{N}(y)$ (recall $E^{\mathcal{A}}_{N}(y)=\mathcal{A}_{N}(\sqrt{y}+E_{N}^{\theta}(y))-\mathcal{A}_{N}(\sqrt{y})$), and using the expressions for $\mathcal{A}^{\pm}(z_{1},..,z_{k})$ from \eqref{e:aplusminus},
we see that the subtraction in the definition of $E_{N}^{\mathcal{A}}$ causes all the terms involving a factor of $\pi^{2}$ to cancel. The terms which remain take one of three forms:
\begin{enumerate}
    \item Terms associated to the $\left(z-\frac{1}{2}\right)\log{(\kappa)}$ term in $\mathcal{A}^{+}$:
    \newline The part of these terms which remains after the cancellations takes the form $\frac{iE_{N}^{\theta}(y)}{2}\log{(\kappa)}$ which is proportional to $E_{N}^{\theta}(y)\log{(N)},$ and therefore, we can apply \cref{l:ethetanbound} to bound this contribution by $CN^{-\chi}(1+\sqrt{r})^{\eta}$ for a small enough $\chi$. 
    \item Terms associated to the $\log{\left(\frac{\Gamma(z)}{\sqrt{2\pi}}\right)}$ term in $\mathcal{A}^{+}:$
    \newline Contributions from these terms after cancellation take the form $\log{(\Gamma(z'_{c}(y)))}-\log{(\Gamma(z_{c}(y)))}$ with $$z_{c}'(y)=c\pm i\frac{\sqrt{y}+E_{N}^{\theta}(y)}{2}, \hspace{10mm} z_{c}(y)=c\pm i \frac{\sqrt{y}}{2}.$$ Therefore, we can apply \cref{l:ethetanbound} (4) to get the appropriate bound. 
    \item Terms associated to the $\left(z-\frac{1}{2}\right)\log{(2)}$ term in $\mathcal{A}^{-}:$
    \newline We note that whenever $N\geq 2$, $\log{(2)}\leq \log{(N)}$ and consequently, since terms of type (3) take the form $E_{N}^{\theta}(y)\log{(2)}$, we can apply \cref{l:ethetanbound} (1) to obtain the desired bound.
\end{enumerate}
\end{proof}
Therefore, using the symmetry of the tangent process \eqref{e:sym}, we conclude that for all $\eta>1$ there exist $N_{0}\in\mathbb{N}$ and  $C,\chi>0$ so that for all $N > N_{0}$ and $y\in [0,4N]$, 
\begin{align} \label{e:conv}
     N^{u+v}\sqrt{1-\frac{y}{4N}}\widehat{\pi}_{t;v,w,u,r}^{(N),c}(y)   = \mathfrak{p}_{t}^{c}(y)\cdot e^{\textsc{Error}_{t}^{(N),c}(y)},
 \end{align} with $|\textsc{Error}_{t}^{(N),c}(y)|  \leq CN^{-\chi}(1+\sqrt{r})^{\eta}.$ This allows us to complete the proof of \cref{l:singlepoint} (1), since
$$\lim_{N\to\infty}N^{u+v}\widehat{\pi}_{t;w,v,u,r}^{(N),c}(y)=C_{b,d,u,v}\cdot \mathfrak{p}_{t}^{c}(y).$$
\end{proof}
We continue to prove \cref{l:singlepoint} (2). In order to prove the lemma, we need to show convergence of the atoms' locations and weights. The convergence result for the locations of the atoms is unchanged from that stated in \cite{CK21}, as the locations of the atoms associated to $\chi\in\{q^{v+\frac{t}{2}},q^{u-\frac{t}{2}}\}$ do not depend on parameters $w$ and $r$, and our choice of time interval excludes all other atoms. 

\begin{lemma}[Lemma 9.3~\cite{CK21}]\label{l:atoms1} Assume $t\in(-2w,2r)$. The location of the atoms of $\widehat{\pi}_{t}^{(N),d}(\cdot)$ associated to $\chi\in\{q^{v+\frac{t}{2}},q^{u-\frac{t}{2}}\}$ converge to those of $\mathfrak{p}_{t}^{d}(\cdot)$. When $v+\frac{t}{2}<0$ the atoms of $\widehat{\pi}_{t}^{(N),d}(\cdot)$ are at $\widehat{y}_{j}^{(N),v}(t)$ for $j\in[0,...,\lfloor -v-\frac{t}{2}\rfloor]$ and $$\lim_{N\to\infty}\widehat{y}_{j}^{(N),v}(t)=x_{j}^{v}(t):=-4\left(v+j+\frac{t}{2}\right)^{2}.$$ Similarly, when $u-\frac{t}{2}<0$ the atoms of $\widehat{\pi}_{t}^{(N),d}(\cdot)$ are at $\widehat{y}_{j}^{(N),u}(t)$ for $j\in[0,...,\lfloor u-\frac{t}{2}\rfloor]$ and $$\lim_{N\to\infty}\widehat{y}_{j}^{(N),u}(t)=x_{j}^{u}(t):=-4\left(u+j+\frac{t}{2}\right)^{2}.$$ 
\end{lemma}

However, it remains to show that their masses converge  to the appropriate atomic masses for the continuous dual Hahn process. The arguments for the atoms arising when $2v+t<0$ are the same as those when the atoms arise from $2u-t<0$, and therefore we will only show the full argument in the former case. 
\begin{proof}[Proof of \cref{l:singlepoint} (2)] We recall the expression for $\widehat{\pi}_{t}^{(N),d}(\widehat{y}_{0}^{(N),v}(t))$,
\begin{align*}
    \widehat{\pi}_{t}^{(N),d}(\widehat{y}_{0}^{(N),v}(t)) & = \frac{(q^{-2v-t}, -q^{w+u}, q^{w+r}, -q^{r+u-t})_{\infty}}{(-q^{w-v}, q^{u-v-t}, -q^{r-v-t}, q^{v+u+w+r})_{\infty}}.    
\end{align*} We can write the logarithm of this expression as
\begin{align*}
    \log{\left(N^{u+v}\widehat{\pi}_{t}^{(N),d}(\widehat{y}_{0}^{(N),v}(t))\right)} & = (u+v)\log{(N)} + \mathcal{A}^{+}_{N}(-2v-t,w+r), + \mathcal{A}^{-}_{N}(w+u,r+u-t)  \\ &  - \mathcal{A}^{+}_{N}(u-v-t,u+v+w+r)-\mathcal{A}^{-}_{N}(w-v,r-v-t) \\ & + E^{+}_{N}(-2v-t,w+r) + E^{-}_{N}(w+u,r+u) \\ &  - E^{+}_{N}(u-v-t,u+v+w+r)-E^{-}_{N}(w-u,r-v-t). 
\end{align*}
The first five terms on the right side simplify to
$$\log{(C_{u,v,w,r})} + \log{(\mathfrak{p}_{t}^{d}(x_{0}^{v}(t)))} = \log{\left(\frac{C_{u,v,w,r}\Gamma(u+v+2)\Gamma(u-v-t)}{\Gamma(-2v-t)}\mathfrak{p}_{t}^{d}(x_{0}^{v}(t))\right)}.$$  
To bound the error terms we simply note that all the arguments of the $E_{N}^{\pm}$ functions are constants, and thus we may apply the bound from \cref{l:ethetanbound} (2) directly. As a result, we see that 
\begin{align}\label{e:discrete0}
N^{u+v}\widehat{\pi}_{t}^{(N),d}(\widehat{y}_{0}^{(N),v}(t)) = C_{u,v,w,r}\cdot \mathfrak{p}_{t}^{d}(x_{0}^{v}(t))\cdot e^{\textsc{Error}_{t}^{d}(x_{0}^{u}(t))},\end{align}
and for all $b\in(0,1)$, there exists a $C>0$ such that  $|\textsc{Error}_{t}^{d}(x_{0}^{u}(t))|\leq CN^{-\frac{b}{2}}.$

The remaining atoms have weights $\widehat{\pi}_{t}^{(N),d}(\widehat{y}_{j}^{(N),v}(t))$,
\begin{align*}
    \widehat{\pi}_{t}^{(N),d}(\widehat{y}_{j}^{(N),v}(t)) & =  \widehat{\pi}_{t}^{(N),d}(\widehat{y}_{0}^{(N),v}(t))\\ & \times \frac{(q^{2v+t},-q^{v+w+t},q^{v+u},-q^{v+r})_{j}(1-q^{2v+2j+t})}{(q,-q^{1+v-w},q^{1+v-u+t},-q^{1+v-r+t})_{j}(1-q^{2v+t})(q^{w+r+u+v-1})^{j}}. 
\end{align*}
showing a similar convergence result for the measure on these atoms reduces to showing a bound on the multiplicative factor. The following lemma follows from equations $(9.58)$ and $(9.59)$ in \cite{CK21}.
\begin{lemma}[p. 58~\cite{CK21}]\label{l:jterms} For all $a\in\mathbb{R}$ if we define functions $E_{N}^{1}(x),E_{N}^{2}(x)$ and $E_{N}^{3}(x)$ such that 
\begin{align*}
    \frac{1-q^{a}}{1-q}=a\cdot e^{E_{N}^{1}(a)}, & & \frac{1+q^{a}}{1+q}=e^{E_{N}^{2}(a)}, & & q^{a}=e^{E^{3}_{N}(a)},
\end{align*} then there exist constants $C_{a}>0$ such that $$|E_{N}^{1}(a)|,|E_{N}^{2}(a)|,|E_{N}^{3}(a)|\leq C_{a}N^{-\frac{1}{2}}.$$
\end{lemma}
Applying this lemma, we can rewrite the multiplicative terms which depend on $j$ as
\begin{align*}
     \frac{(v+j+\frac{t}{2})[2v+t,v+u]_{j}}{(v+\frac{t}{2})j![1+v-u+t]_{j}}e^{E_{N}^{M_{j}}},
\end{align*} and there exist $C\in\mathbb{R}_{>0}$ such that $|E_{N}^{M_{j}}|\leq CN^{-\frac{1}{2}}$. Therefore, applying the definition of $\mathfrak{p}_{t}^{d}(x_{j}^{v}(t))$ from \cref{d:cdhproc}, we see that 
\begin{align}\label{e:discretej}N^{u+v}\widehat{\pi}_{t}^{v}(y_{j}^{(N),d})=C_{u,v,w,r}\cdot\mathfrak{p}_{t}^{d}(x_{j}^{v}(t)).\end{align}
The argument for the atoms arising from $u-\frac{t}{2}$ is identical, and will therefore not be repeated. 

As a result of \eqref{e:discrete0} and \eqref{e:discretej} and the analogous results for the atoms arising from $u-\frac{t}{2}$, we see that 
\begin{align*}
    \lim_{N\to\infty}N^{u+v}\widehat{\pi}_{t}^{v}(y_{j}^{(N),d}) & =C_{u,v,w,r}\cdot\mathfrak{p}_{t}^{d}(x_{j}^{v}(t)), \\ \lim_{N\to\infty}N^{u+v}\widehat{\pi}_{t}^{u}(y_{j}^{(N),d}) & =C_{u,v,w,r}\cdot\mathfrak{p}_{t}^{d}(x_{j}^{u}(t)), 
\end{align*} which completes the proof of \cref{l:singlepoint} (2). \end{proof}

\subsection{Proof of \cref{l:transprob}}\label{s:transconv}

\begin{proof} [Proof of \cref{l:transprob} (1)]
    We begin by recalling the expression for $\widehat{\pi}_{s,t;v,w,u,r}^{(N),c,c}(m,r)$: 
    \begin{align}\label{e:cctransform}\widehat{\pi}_{s,t;v,w,u,r}^{(N),c,c}(m,r)=\frac{1}{2N}\pi_{q^{t},q^{s}}^{(N),c,c}\left(1-\frac{r}{2N}, 1-\frac{m}{2N}\right)\cdot\frac{\widehat{\pi}_{t;v,w,u,r}^{(N),c}(r)}{\widehat{\pi}_{s;v,w,u,r}^{(N),c}(m)}.
    \end{align} 

Therefore, we need to understand the asymptotics of the following Askey-Wilson measure expression. 
\begin{align*}
& \sqrt{1-\frac{m}{4N}}\frac{1}{2N}\pi_{q^{t},q^{s}}\left(1-\frac{r}{2N},1-\frac{m}{2N}\right)   \\ & =
\frac{\left(q,-q^{v+w+s}, q^{t-s}\right)_{\infty}\bigg|\left(q^{v+\frac{t}{2}+i\frac{\sqrt{N}\theta_{r}}{2}}, -q^{w+\frac{t}{2}+i\frac{\sqrt{N}\theta_{r}}{2}},q^{i\frac{\sqrt{N}\theta_{m}}{2}}\right)\bigg|^{2}}{4\pi N\sqrt{\frac{m}{N}} \left(-q^{v+w+t}\right)_{\infty}\bigg|\left(q^{v+\frac{s}{2}+i\frac{\sqrt{N}\theta_{m}}{2}}, -q^{w+\frac{s}{2}+i\frac{\sqrt{N}\theta_{m}}{2}}, q^{\frac{t-s}{2} + i\frac{\sqrt{N}(\theta_{r}+\theta_{m})}{2}},q^{\frac{t-s}{2} + i\frac{\sqrt{N}(-\theta_{r}+\theta_{m})}{2}}\right)_{\infty}\bigg|^{2}}.
\end{align*}

With a slight modification of the notation we used in \cref{s:sinconv}, defined below, we can rewrite the logarithm of this equation.

\begin{equation}\label{e:logcc}
\begin{aligned}
    & \log{\left(\sqrt{1-\frac{m}{4N}}\frac{1}{2N}\pi_{q^{t},q^{s}}^{(N),c,c}\left(1-\frac{r}{2N}, 1-\frac{m}{2N}\right)\right)}  \\ & =  -\log{\left(4\pi N\sqrt{\frac{m}{N}}\right)} + \mathcal{A}_{N}(\sqrt{m},\sqrt{r})   + E_{N}^{\mathcal{A}}(m,r) + E_{N}^{\mathcal{E}}(m,r) .
\end{aligned}
\end{equation}
In the expression above, 
\begin{align*} \mathcal{A}_{N}(\sqrt{m},\sqrt{r})& := \mathcal{A}_{N}^{+}\left(1,t-s, v+\frac{t}{2} + i\frac{\sqrt{r}}{2}, v+\frac{t}{2} - i\frac{\sqrt{r}}{2}, i\sqrt{m}, -i\sqrt{m}\right) \\ & + \mathcal{A}_{N}^{-}\left(v+w+s,w+\frac{t}{2}+i\frac{\sqrt{r}}{2},w+\frac{t}{2}-i\frac{\sqrt{r}}{2}\right) \\ &  -\mathcal{A}^{+}_{N}\left(v+\frac{s}{2} + i\frac{\sqrt{m}}{2},v+\frac{s}{2} - i\frac{\sqrt{m}}{2}, \frac{t-s}{2} + i \frac{\sqrt{r}+\sqrt{m}}{2}\right) \\ & -\mathcal{A}^{+}_{N}\left(\frac{t-s}{2} - i \frac{\sqrt{r}+\sqrt{m}}{2},\frac{t-s}{2} + i \frac{\sqrt{r}-\sqrt{m}}{2},\frac{t-s}{2} - i \frac{\sqrt{r}-\sqrt{m}}{2}\right) \\ & -\mathcal{A}^{-}_{N}\left(v+w+t, w+\frac{s}{2} + i\frac{\sqrt{m}}{2}, w+\frac{s}{2} - i\frac{\sqrt{m}}{2}\right), 
\\ E_{N}^{\mathcal{A}}(m,r) & :=\mathcal{A}_{N}(\sqrt{m}+E_{N}^{\theta}(m),\sqrt{r}+E_{N}^{\theta}(r))-\mathcal{A}_{N}(\sqrt{m},\sqrt{r}), 
\\ E^{\mathcal{E}}_{N}(m,r) & :=\mathcal{E}_{N}(\sqrt{m}+E^{\theta}_{N}(m),\sqrt{r}+E^{\theta}_{N}(r)).\end{align*}

As before, $E_{N}^{\theta}(a)$ is defined by $\sqrt{N}\theta_{a}^{(N)}=\sqrt{a}+E_{N}^{\theta}(a).$ We can simplify the first two terms of \eqref{e:logcc} by expanding using the definitions of $\mathcal{A}^{\pm}$, \eqref{e:aplusminus}, and then canceling like terms.

\begin{align*}
    & -\log{\left(4\pi N\sqrt{\frac{r}{N}}\right)}+\mathcal{A}_{N}(\sqrt{m},\sqrt{r}) \\ &  =\log{\left(\frac{\bigg|\Gamma\left(v+\frac{s}{2} + i\frac{\sqrt{m}}{2}\right)\Gamma\left(\frac{t-s}{2} + i \frac{\sqrt{r}+\sqrt{m}}{2}\right)\Gamma\left(\frac{t-s}{2} + i \frac{\sqrt{r}-\sqrt{m}}{2}\right)\bigg|^{2}}{8\pi \sqrt{r}\Gamma(t-s)\bigg|\Gamma\left(v+\frac{t + i\sqrt{r}}{2}\right)\Gamma\left( \frac{i\sqrt{m}}{2}\right)\bigg|^{2}}\right)}.
\end{align*}

We also simplify the expression for the ratio of the marginal distributions which arises in \eqref{e:cctransform}. 

$$\log{\left(\frac{\sqrt{1-\frac{r}{4N}}\widehat{\pi}_{t;v,w,u,r}^{(N),c}(r)}{\sqrt{1-\frac{m}{4N}}\widehat{\pi}_{s;v,w,u,r}^{(N),c}(m)}\right)}=\log{\left(\frac{\mathfrak{p}_{t}^{c}(r)}{\mathfrak{p}_{s}^{c}(m)}\right)}+ \textsc{Error}_{t}^{(N),c}(r)-\textsc{Error}_{s}^{(N),c}(m),$$

using these expressions, and the definitions of the continuous dual Hahn process, \cref{d:cdhproc}, we can rewrite \eqref{e:cctransform} as

\begin{align*}
    \log{\left(\sqrt{1-\frac{r}{4N}}\pi_{s,t}^{(N),c,c}(m,r)\right)} &
 =  E_{N}^{\mathcal{A}}(m,r) + E_{N}^{\mathcal{E}}(m,r) + \textsc{Error}_{t}^{(N),c}(r) - \textsc{Error}_{s}^{(N),c}(m)
\end{align*}
\begin{align*}+\log{\left(\frac{\bigg|\Gamma\left(v+\frac{s}{2} + i\frac{\sqrt{m}}{2}\right)\Gamma\left(\frac{t-s}{2} + i \frac{\sqrt{r}+\sqrt{m}}{2}\right)\Gamma\left(\frac{t-s}{2} + i\frac{\sqrt{r}-\sqrt{m}}{2}\right) \bigg|^{2}\mathfrak{p}_{t}^{c}(r)}{8\pi \sqrt{r}\Gamma(t-s)\bigg|\Gamma\left(v+\frac{t}{2} + i\frac{\sqrt{r}}{2}\right)\Gamma\left( i\frac{\sqrt{m}}{2}\right)\bigg|^{2}\mathfrak{p}_{s}^{c}(m)}\right)} 
\end{align*}
\begin{align*} = \log{\left(\mathfrak{p}_{s,t}^{c,c}(m,r)\right)} + E_{N}^{\mathcal{A}}(m,r) + E_{N}^{\mathcal{E}}(m,r) +\textsc{Error}_{t}^{(N),c}(r) - \textsc{Error}_{s}^{(N),c}(m).
\end{align*}
The last step is to bound the error terms. The bounds on $\textsc{Error}_{t}^{(N),c}(r)$ and $\textsc{Error}_{s}^{(N),c}(m)$ were already given in \cref{s:sinconv}. The remaining error bound is stated in the following lemma.
\begin{lemma}
    For all $\eta>1$ there exists $N_{0}\in\mathbb{Z}_{\geq 1}$ and $C,\chi>0$ such that for all $N>N_{0}$ and $m,r\in\widehat{S}^{(N),c}=[0,4N]$, $$|E_{N}^{\mathcal{A}}(m,r)|,|E_{N}^{\mathcal{E}}(m,r)|\leq CN^{-\chi}(1+\sqrt{m})^{\eta}+CN^{-\chi}(1+\sqrt{r})^{\eta}.$$
\end{lemma}
The proof of this lemma is the same as the proof of \cref{p:ebounds}, the only change being the slight difference in definition of the error functions, so we will not repeat it.
\end{proof}

\begin{proof}[Proof of \cref{l:transprob} (2)]
We rewrite $\widehat{\pi}_{s,t;v,w,u,r}^{(N),d,c}\left(\widehat{y}_{j}^{(N),\chi}(s),r\right)$ through the same time reversal transformation:

\begin{align}\label{e:dcvtimerev}
\widehat{\pi}_{s,t;v,w,u,r}^{(N),d,c}(\widehat{y}_{j}^{(N),\chi}(s),r)=\pi_{q^{t},q^{s}}^{(N),c,d}\left(1-\frac{r}{2N},y_{j}^{(N),\chi}(s)\right)\cdot\frac{\widehat{\pi}_{t;v,w,u,r}^{(N),c}(r)}{\widehat{\pi}_{s}^{(N),d}(\widehat{y}_{j}^{(N),\chi}(s))}.
\end{align}

We analyze the expression
$$\pi_{q^{t},q^{s}}^{(N)}\left(1-\frac{r}{2N},V\right)=AW(V; q^{v+\frac{s}{2}}, -q^{w+\frac{s}{2}}, q^{\frac{t-s}{2} + i\frac{\sqrt{N}\theta_{r}}{2}}, q^{\frac{t-s}{2} - i\frac{\sqrt{N}\theta_{r}}{2}}).$$

Atoms may arise when $2v+s<0$, with associated weights
    \begin{align}\label{e:dcv}
\pi_{q^{t},q^{s}}^{(N),c,d}\left(1-\frac{r}{2N}, y_{0}^{(N),v}(s)\right)
& = \frac{(q^{-2v-s}, -q^{w+\frac{t}{2}+i\frac{\sqrt{N}\theta_{r}}{2}}, -q^{w+\frac{t}{2} - i\frac{\sqrt{N}\theta_{r}}{2}}, q^{t-s})_{\infty}}{(-q^{w-v}, q^{-v+\frac{t}{2}-s+i\frac{\sqrt{N}\theta_{r}}{2}}, q^{-v+\frac{t}{2}-s-i\frac{\sqrt{N}\theta_{r}}{2}}, -q^{v+t+w})_{\infty}},
\\ \pi_{q^{t},q^{s}}^{(N),c,d}\left(1-\frac{r}{2N}, y_{j}^{(N),v}(s)\right)
& = \pi_{q^{t},q^{s}}^{(N),c,d}\left(1-\frac{r}{2N}, y_{0}^{(N),v}(s)\right) \times M_{j},
\end{align}
\begin{align}\label{e:dcvm}
    M_{j} & = \frac{(q^{2v+s}, -q^{v+w+s}, q^{v+\frac{t}{2} + i\frac{\sqrt{N}\theta_{r}}{2}}, - q^{v+\frac{t}{2}-i\frac{\sqrt{N}\theta_{r}}{2}})_{j}(1-q^{2v+2j+s})}{(q, -q^{1+v-w},q^{1+v+s -\frac{t}{2} -i\frac{\sqrt{N}\theta_{r}}{2}},q^{1+v+s -\frac{t}{2} +i\frac{\sqrt{N}\theta_{r}}{2}})_{j}(1-q^{2v+s})(-q^{v+w+t-1})^{j}}.
\end{align}

To analyze the weights which arise from the atoms, we need to ensure that there are no instances of division by zero. Any division by zero in \eqref{e:dcv} arises when $-v+\frac{t}{2}-s\in\mathbb{Z}_{\leq0}$ or, in \eqref{e:dcvm}, when $1+v+s-\frac{t}{2}+j-1=0$ for $j\in[[1,\lfloor -v-\frac{s}{2}\rfloor]]$, or when $2v+s=0$. The final case can never happen, since we assumed, in order to produce atoms, that $v+\frac{s}{2}<0$. Furthermore, since $s<t$, the other two scenarios are also ruled out. We see that $$\log{\left(\pi_{q^{t},q^{s}}^{(N),c,d}\left(1-\frac{r}{2N}, y_{0}^{(N),v}(s)\right)\right)}=\mathcal{A}_{N}(\sqrt{r})+E_{N}^{\mathcal{A}}(r)+E_{N}^{\mathcal{E}}(r).$$
Where, in this equation, we define $\mathcal{A}_{N}(a)$ for $a\in\mathbb{R}_{\geq 0}$ by
\begin{equation}\label{e:es}
\begin{aligned}
    \mathcal{A}_{N}(\sqrt{r}) & := \mathcal{A}_{N}^{+} \left(-2v-s, t-s\right) + \mathcal{A}_{N}^{-}\left(b+\frac{t}{2}+i\frac{\sqrt{r}}{2}, b+\frac{t}{2}-i\frac{\sqrt{r}}{2}\right) \\ & -\mathcal{A}_{N}^{+}\left(-v+\frac{t}{2}-s + i\frac{\sqrt{r}}{2}, -v+\frac{t}{2} - s-i\frac{\sqrt{r}}{2}\right) - \mathcal{A}_{N}^{-}\left(b-v, v+t+b\right)  \\ & =\log{\left(\frac{\bigg|\Gamma\left(\frac{t}{2} -s -v+\frac{i\sqrt{r}}{2}\right)\bigg|^{2}}{\Gamma(-2v-s)\Gamma(t-s)}\right)},  \\
    E_{N}^{\mathcal{A}}(r) & :=\mathcal{A}_{N}(\sqrt{r}+E^{\theta}(r))-\mathcal{A}_{N}(\sqrt{r}), \\ E_{N}^{\mathcal{E}}(r) & :=\mathcal{E}_{N}(\sqrt{r}+E_{N}^{\theta}(r)).
\end{aligned}
\end{equation}

We also need to address the ratio of $\widehat{\pi}_{t}^{(N),c}(r)$ and $\widehat{\pi}_{s}^{(N),d}(y_{j}^{(N),v}(s))$ which arose in \eqref{e:dcvtimerev}. 
\begin{align*}
    \log{\left(\frac{\sqrt{1-\frac{r}{4N}}\widehat{\pi}_{t;v,w,u,r}^{(N),c}(r)}{\widehat{\pi}_{s}^{(N),d}(\widehat{y}_{j}^{(N),v}(s))}\right)} = \log{\left(\frac{\mathfrak{p}_{t}^{c}(y)}{\mathfrak{p}_{s}^{d}(x_{j}^{v})}\right)} +\textsc{Error}_{t}^{(N),c}(r) - \textsc{Error}_{s}^{(N),d}(x_{0}^{v}(s)).
\end{align*} Consequently, applying the definitions of the continuous dual Hahn process measures from \cref{d:cdhproc}, we find
\begin{equation}\label{e:errorterms}
\begin{aligned}
    \log{\left(\sqrt{1-\frac{r}{4N}}\widehat{\pi}_{s,t;v,w,u,r}^{(N),d,c}(\widehat{y}_{0}^{(N),v}(s),r)\right)} & = \log{\left(\frac{|\Gamma(\frac{t}{2}-s-v+\frac{i\sqrt{r}}{2})|^{2}\mathfrak{p}_{t}^{c}(r)}{\Gamma(-2v-s)\Gamma(t-s)\mathfrak{p}_{s}^{d}(x_{j}^{v})}\right)}   + E_{N}^{\mathcal{A}}(r) \\ &  + E_{N}^{\mathcal{E}}(r)  + \textsc{Error}_{t}^{(N),c}(r) -\textsc{Error}_{s}^{(N),d}(x_{0}^{v}(s))  \\ & = \log{\left(\mathfrak{p}_{s,t}^{d,c}(x_{0}^{v}(s),r)\right)}  + E_{N}^{\mathcal{A}}(r)+E_{N}^{\mathcal{E}}(r) \\ & + \textsc{Error}_{t}^{(N),c}(r)- \textsc{Error}_{s}^{(N),d}(x_{0}^{v}(s)) .
\end{aligned}
\end{equation}
To finish the argument, we need the following error bounds. 
\begin{lemma}\label{l:moreerrorbounds}
    For $E_{N}^{\mathcal{A}}(r)$ and $E_{N}^{\mathcal{E}}(r)$ as defined by \eqref{e:es}, for all $\eta>1$ there exists $N_{0}\in\mathbb{Z}_{\geq 1}$ and $C,\chi>0$ such that for all $N>N_{0}$ and $r\in\widehat{S}^{(N),c}=[0,4N]$, 
    $$|E_{N}^{\mathcal{A}}(r)|,|E_{N}^{\mathcal{E}}(r)|\leq CN^{-\chi}(1+\sqrt{r})^{\eta}.$$
\end{lemma}
The proof of this lemma is unchanged from the earlier bound on the error terms in \cref{p:ebounds}, the only change being the modified definitions of $E_{N}^{\mathcal{A}}(r)$ and $E_{N}^{\mathcal{E}}(r)$, and therefore we will not repeat the proof. Combining the bounds from \cref{l:moreerrorbounds} and \cref{l:ethetanbound}, we see that for all $\eta>1$, there exist $C,\chi>0$ such that the error in \eqref{e:errorterms} is bounded by $CN^{-\chi}(1+\sqrt{r})^{\eta}$. This bound completes the argument in the case that $j=0$, and all that remains is to investigate the terms in $M_{j}$, defined in \eqref{e:dcvm}. 
\begin{lemma}[p. 64~\cite{CK21}]\label{l:multiplicative} The multiplicative term $M_{j}$ can be written as 
\begin{align}
    M_{j}=\frac{[2v+s,v+\frac{t}{2}+i\frac{\sqrt{r}}{2},v+\frac{t}{2}-i\frac{\sqrt{r}}{2}]_{j}(2v+2j+s)}{[1, 1+v+s-\frac{t}{2}+i\frac{\sqrt{r}}{2},1+v+s-\frac{t}{2}-i\frac{\sqrt{r}}{2}]_{j}(2v+s)}(-1)^{j}e^{E_{N}^{M_{j}}(r)}.
\end{align} Where, for all $\eta>1$ there exists $N_{0}\in\mathbb{Z}_{\geq 1}$ and $C,\chi>0$ such that for all $N>N_{0}$ and $r\in\widehat{S}_{t}^{(N),c}=[0,4N]$, $$|E_{N}^{M_{j}}(r)|\leq CN^{-\chi}(1+\sqrt{r})^{\eta}.$$
\end{lemma}

Applying this lemma, we see that $$\log{\left(\sqrt{1-\frac{r}{4N}}\widehat{\pi}_{s,t;v,w,u,r}^{(N),d,c}(\widehat{y}_{j}^{(N),v}(s),r)\right)} = \log{\left(\mathfrak{p}_{s,t}^{d,c}(x_{j}^{v},r)\right)}+\overline{\textsc{Error}}_{s,t}^{(N),d,c}(x_{j}^{v}(s),r),$$ with $\overline{\textsc{Error}}_{s,t}^{(N),d,c}(x_{j}^{v}(s),r):=\textsc{Error}_{s,t}^{(N),d,c}(x_{j}^{v}(s),r)+E_{N}^{M_{j}}$. By the above calculations, for all $\eta>1$, there exist $N_{0}\in\mathbb{Z}_{\geq 1}$ and $C,\chi>0$ such that for all $N>N_{0}$, $r\in \widehat{S}^{(N),c}_{t}$, $$\overline{\textsc{Error}}_{s,t}^{(N),d,c}(x_{j}^{v}(s),r) \leq CN^{-\chi}(1+\sqrt{r})^{\eta}.$$
\end{proof}
Finally, we can give the proof of \cref{l:transprob} (3). 
\begin{proof}[Proof of \cref{l:transprob} (3)] Here $u-\frac{s}{2}<0$ and consequently, $v+\frac{s}{2}\geq 0$. We note that the Askey-Wilson measure associated to this case will be 

$$\pi_{q^{t},q^{s}}^{(N)}\left(1-\frac{r}{2N},V\right)=AW(V; q^{v+\frac{s}{2}}, -q^{w+\frac{s}{2}}, q^{\frac{t-s}{2} + i\frac{\sqrt{N}\theta_{r}}{2}}, q^{\frac{t-s}{2} - i\frac{\sqrt{N}\theta_{r}}{2}}).$$ 
The only possibility is that we have atoms arising from $w+\frac{t}{2}<0$, but we have restricted our domain in time so that these atoms do not arise. 
\end{proof}

\subsection{Proof of \cref{l:ddconv}}\label{s:transdconv} 

We begin with the proof of part (1).
\begin{proof}[Proof of \cref{l:ddconv} (1)] We know that the time-reversed transition measure can be expressed via the following transformation on the standard Askey-Wilson measure. 
$$\widehat{\pi}_{s,t}^{(N),d,d}(x,y)=\pi_{q^{t},q^{s}}^{(N),d,d}\left(1-\frac{y}{2N},1-\frac{x}{2N}\right)\cdot\frac{\widehat{\pi}_{t}^{(N),d}(y)}{\widehat{\pi}_{s}^{(N),d}(x)}.$$
The relevant Askey-Wilson measure is
\begin{align*}
    \pi_{q^{t},q^{s}}^{(N)}\left(y_{j}^{(N),v}(t), V\right) & = AW\left(V;q^{v+\frac{s}{2}},-q^{w+\frac{s}{2}}, q^{\frac{t}{2}+v+j}, q^{\frac{t}{2}-s-v-j}\right).
\end{align*}
This measure will only have atoms arising from $v+\frac{s}{2}<0$ and $w+\frac{s}{2}<0$, however, we restrict the domain in $t$ so that the latter case never occurs. What remains are atoms arising from $v+\frac{s}{2}<0$, which have weights
\begin{align*}
  \pi_{q^{t},q^{s}}^{(N)}\left(y_{j}^{(N),v}(t),y_{0}^{(N),v}(s)\right) 
    & = \frac{(q^{-2v-s}, -q^{w+t +j+v}, -q^{w -j-v}, q^{t-s})_{\infty}}{(-q^{w-v}, q^{t -s + j}, q^{-2v -s -j}, -q^{v+w+t})_{\infty}} ,
    \\
    \pi_{q^{t},q^{s}}^{(N)}\left(y_{j}^{(N),v},y_{k}^{(N),v}\right) & = \pi_{q^{t},q^{s}}^{(N)}\left(y_{j}^{(N),v},y_{0}^{(N),v}\right)\times M_{k}, \\ M_{k} & = \frac{(q^{2v+s},-q^{v+w+s}, q^{2v+t+ j}, q^{-j})_{k}(1-q^{2v+2k+s})}{(q, -q^{1+v-w}, q^{1+s - t-j},q^{1+2v+s+j})_{k}(1-q^{2v+s})(-q^{v+w+t-1})^{k}}.
\end{align*}

Taking the logarithm, we write
\begin{align*}
    \log{\left(\pi_{q^{t},q^{s}}^{(N)}\left(y_{j}^{(N),v}(t),y_{0}^{(N),v}(s)\right)\right)} & = \mathcal{A}_{N}(j) + E_{N}^{\mathcal{A}}(j) + E_{N}^{\mathcal{E}}(j),
\end{align*}
where $\mathcal{A}_{N}(j)$ is defined
\begin{align*}
    \mathcal{A}_{N}(j) & := \mathcal{A}^{+}_{N}\left(-2v-s, t-s\right)+ \mathcal{A}^{-}\left(w+t+j+v,w-j-v\right) \\ & - \mathcal{A}^{+}\left(t-s+j,-2v-s-j\right) - \mathcal{A}^{-}\left(w-v,v+w+t\right)
    \\ & =   \log{\left(\frac{\Gamma(t-s+j)\Gamma(-2v-s-j)}{\Gamma(-2v-s)\Gamma(t-s)}\right)}. 
\end{align*}
As before, albeit with modified error functions $E_{N}^{\mathcal{A}}(j)$ and $E_{N}^{\mathcal{E}}(j)$, we have the following lemma. 
\begin{lemma}
    There exists $N_{0}\in\mathbb{Z}_{\geq 1}$ and $C,\chi>0$ such that for all $N>N_{0}$ and $j\in[[0,\lfloor -u-\frac{t}{2}\rfloor]]$, $E_{N}^{\mathcal{A}}(j)$ and $E_{N}^{\mathcal{E}}(j)$ are bounded by $$|E_{N}^{\mathcal{A}}(j)|,|E_{N}^{\mathcal{E}}(j)|\leq CN^{-\chi}.$$
\end{lemma}

Using the lemma, we find that 
$$\pi_{q^{t},q^{s}}^{(N),d,d}\left(y_{j}^{(N),v}(t),y_{0}^{(N),v}(s)\right)=\frac{\Gamma(t-s+j)\Gamma(-2v-s-j)}{\Gamma(-2v-s)\Gamma(t-s)}e^{\textsc{Error}_{t,s}^{(N),d,d}(y_{j}^{(N),v}(t),y_{0}^{(N),v}(s))},$$
and there exist $C,\chi>0$ and $N_{0}\in\mathbb{Z}_{\geq 1}$ so that for all $N>N_{0}$, $$|\textsc{Error}_{s,t}^{(N),d,d}(y_{j}^{(N),v}(t),y_{0}^{(N),v}(s))|\leq CN^{-\chi}.$$

Now, we just need to deal with the $M_{k}$ terms. We can apply a version of \cref{l:multiplicative} to conclude that  
$$M_{k} =\frac{[2v+s, 2v+t+j, -j]_{k}(2v+2k+s)}{[1, 1+s-t-j,1+2v+s+j]_{k}(2v+s)}(-1)^{k}e^{\overline{\textsc{Error}}_{t,s}^{(N),d,d}(y_{j}^{(N),v}(t),y_{k}^{(N),v}(s))}.$$ and there exist $C,\chi>0$ and $N_{0}\in\mathbb{Z}_{\geq 1}$ such that for all $N>N_{0}$, 
$$|\overline{\textsc{Error}}_{t,s}^{(N),d,d}(y_{j}^{(N),v}(t),y_{k}^{(N),v}(s))|\leq CN^{-\chi}.$$

We can also deal with the limit of the ration of $\widehat{\pi}_{t}^{(N),d}(y)$ and $\widehat{\pi}_{s}^{(N),d}(x)$ using our marginal computations. 
\begin{align*}
    \lim_{N\to\infty}\frac{\widehat{\pi}_{t}^{(N),d}(\widehat{y}_{j}^{(N),v})}{\widehat{\pi}_{s}^{(N),d}(\widehat{y}_{j}^{(N),v})} 
 & = \lim_{N\to\infty}\frac{N^{u+v}\widehat{\pi}_{t}^{(N),d}(\widehat{y}_{j}^{(N),v}(t))}{N^{u+v}\widehat{\pi}_{s}^{(N),d}(\widehat{y}_{k}^{(N),v}(s))}  = \frac{\mathfrak{p}_{t}^{d}(x_{j}^{v}(t))}{\mathfrak{p}_{s}^{d}(x_{j}^{v}(s))}.
\end{align*}

Consequently, we can prove \cref{l:ddconv} (1). 
$$\lim_{N\to\infty}\widehat{\pi}_{s,t}^{(N),d,d}(\widehat{y}_{j}^{(N),v},\widehat{y}_{k}^{(N),v}) = \mathfrak{p}_{s,t}(x_{k}^{v}(s),x_{j}^{v}(t)) .$$ 
\end{proof}

We now turn to \cref{l:ddconv} (2). 

\begin{proof}[Proof of \cref{l:ddconv} (2)]
We begin by considering the transition measures, 
\begin{align*}
    \widehat{\pi}_{s,t}^{(N),d,d}\left(\widehat{y}_{j}^{(N),u}(s),\widehat{y}_{k}^{(N),v}(t)\right) &  = \pi_{q^{t},q^{s}}^{(N),d,d}\left(y_{k}^{(N),v}(t),y_{j}^{(N),u}(s)\right)\frac{\pi_{q^{t}}^{(N),d}\left(y_{k}^{(N),v}(t)\right)}{\pi_{q^{s}}^{(N),d}\left(y_{j}^{(N),u}(s)\right)}, \\ 
    \pi_{q^{t},q^{s}}^{(N),d,c}\left(y_{k}^{(N),v}(t), V\right) & = AW\left(V; q^{v+\frac{s}{2}}, -q^{w+\frac{s}{2}}, q^{\frac{t-s}{2} + \left(k+v+\frac{t}{2}\right)}, q^{\frac{t-s}{2} -\left(k+v+\frac{t}{2}\right)}\right),
\end{align*} which has no $u$ atoms. Likewise,

\begin{align*} \widehat{\pi}_{s,t}^{(N),d,d}\left(\widehat{y}_{j}^{(N),u}(s),\widehat{y}_{k}^{(N),u}(t)\right) & = \pi_{q^{t},q^{s}}^{(N),d,d}\left(y_{k}^{(N),u}(t),y_{j}^{(N),u}(s)\right)\frac{\pi_{q^{t}}^{(N),d}\left(y_{k}^{(N),u}(t)\right)}{\pi_{q^{s}}^{(N),d}\left(y_{j}^{(N),u}(s)\right)} \\ \pi_{q^{t},q^{s}}^{(N),d,d}\left(y_{k}^{(N),u}(t),V\right) & = AW\left(V;q^{v+\frac{s}{2}}, -q^{w+\frac{s}{2}}, q^{\frac{t-s}{2}+u+k-\frac{t}{2} },q^{\frac{t0s}{2}-(u+k-\frac{t}{2})}\right),
\end{align*} which no atoms associated to the first and second arguments, due to the constraint $u+v>0$ and $u-s/2<0$ as well as the restriction of the domain to avoid atoms arising from $w+s/2<0.$ Only the third argument can attain atoms, which arise whenever $-s/2 +u +k<0 $, in which case the term will contribute an atom $y^{(N),u}_{k+i}(s)$ for all $i\in [[0,\lfloor -u-k+\frac{s}{2}\rfloor]]$. The masses at these atoms are given by 
\begin{align*}
    \pi_{q^{t},q^{s}}^{(N),d,d}\left(y_{k}^{(N),u}(t),y_{k}^{(N),u}(s)\right) & = \frac{(q^{-2u-2k+s},-q^{v+1+s},q^{v-u-k+t},-q^{1-u-k+t})_{\infty}}{(q^{v-u-k+s},-q^{1-u-k+s},q^{-2u-2k+t},-q^{v+t+1})_{\infty}},
    \\ \pi_{q^{t},q^{s}}^{(N),d,d}\left(y_{k}^{(N),u}(t),y_{k+i}^{(N),u}(s)\right) & = \pi_{q^{t},q^{s}}^{(N),d,d}\left(y_{k}^{(N),u}(t),y_{k}^{(N),u}(s)\right) \cdot M_{i}
    \\ M_{i} & := \frac{(q^{2u+2k-s},q^{v+u+k},-q^{u+1+k},q^{t-s})_{i}(1-q^{2u+2k+2i-s})}{(q,q^{u-v+1+k-s},-q^{u+k-s},q^{2u+2k-t+1})_{i}(1-q^{2u+2k-s})}\cdot (-q^{-v-2t+s})^{i}
\end{align*} 
Under the same analysis as in preceding discussions, we see that 
\begin{align*}
    \lim_{N\to\infty}\pi_{q^{t},q^{s}}^{(N),d,d}\left(y_{k}^{(N),u}(t),y_{k}^{(N),u}(s)\right) & = \frac{\Gamma(v-u-k+s)\Gamma(-2u-2k+t)}{\Gamma(-2u-2k+s)\Gamma(v-u-k+t)}
    \\ \lim_{N\to\infty} M_{i} & = \frac{[2u+2k-s,v+u+k,t-s]_{i}(2u+2k+2i-s)}{[1,u-v+1+k-s,2u+2k-t+1]_{i}(2u+2k-s)}(-1)^{i}.
\end{align*}
Putting this together, and replicating the asymptotic bounds from previous computations finishes the proof of \cref{l:ddconv} (2). 
\end{proof}

Finally, we prove \cref{l:ddconv} (3).
\begin{proof}[Proof of \cref{l:ddconv} (3)]
We write
\begin{align*}
    \frac{1}{2N}\pi_{q^{t},q^{s}}^{(N)}\left(y_{j}^{(N),v}(t), 1-\frac{r}{2N}\right) & = \frac{1}{4\pi N\sqrt{\frac{r}{N}}\sqrt{1-\frac{r}{4N}}}
\end{align*}
\begin{align*}\times 
\frac{(q, -q^{v+w+s}, q^{t+j+2v},q^{-j},-q^{t+w+j+v},-q^{w-j-v},q^{t-s})_{\infty}|(q^{i\sqrt{N}\theta_{r}})|^{2}}{(-q^{v+w+t})_{\infty}\bigg|\left(q^{v+\frac{s}{2}+i\frac{\sqrt{N}\theta_{r}}{2}}, -q^{w+\frac{s}{2}+i\frac{\sqrt{N}\theta_{r}}{2}}, q^{\frac{t-s}{2} + (j+v+\frac{t}{2})+i\frac{\sqrt{N}\theta_{r}}{2}}, q^{\frac{t-s}{2} -(j+v+\frac{t}{2})+i\frac{\sqrt{N}\theta_{r}}{2}}\right)_{\infty}\bigg|^{2}}.
\end{align*}

The $(q^{-j})_{\infty}$ term in the numerator sends the entire expression to $0$. We can consider instead the transition from the atoms associated to $u$. 
\begin{align*}
    \frac{1}{2N}\pi_{q^{t},q^{s}}^{(N)}\left(y_{j}^{(N),u}(t), 1-\frac{r}{2N}\right) & = \frac{1}{4\pi N\sqrt{\frac{r}{N}}\sqrt{1-\frac{r}{4N}}} 
\end{align*}
\begin{align*}\times 
 \frac{(q, -q^{v+w+s}, q^{j+v+u},q^{v-j-u+t},-q^{w+j+u},-q^{w-j-u+t},q^{t-s})_{\infty}|(q^{i\sqrt{N}\theta_{r}})|^{2}}{(-q^{v+w+t})_{\infty}\bigg|\left(q^{v+\frac{s}{2}+i\frac{\sqrt{N}\theta_{r}}{2}}, -q^{w+\frac{s}{2}+i\frac{\sqrt{N}\theta_{r}}{2}}, q^{\frac{t-s}{2} + (j+u-\frac{t}{2})+i\frac{\sqrt{N}\theta_{r}}{2}}, q^{\frac{t-s}{2} -(j+u-\frac{t}{2})+i\frac{\sqrt{N}\theta_{r}}{2}}\right)_{\infty}\bigg|^{2}}.
\end{align*} With the same asymptotic analysis as in prior sections, we find that 
\begin{align*} & \sqrt{1-\frac{r}{4N}}\frac{1}{2N}\pi_{q^{t},q^{s}}^{(N)}\left(y_{j}^{(N),u}(t),1-\frac{r}{2N}\right)
\\ & =\frac{\big|\Gamma\left(v+\frac{s}{2}+i\frac{\sqrt{r}}{2}\right)\Gamma\left(j+u-\frac{s}{2}+i\frac{\sqrt{r}}{2}\right)\Gamma\left(t-\frac{s}{2}-j-u+i\frac{\sqrt{r}}{2}\right)\big|^{2}}{8\pi\sqrt{r}\Gamma\left(j+v+u\right)\Gamma\left(v-j-u+t\right)\Gamma\left(t-s\right)|\Gamma(i\sqrt{r})|^{2}} \cdot e^{\textsc{Error}_{t,s}^{(N)}(y^{(N),u}_{j}(t))},
\end{align*} with the error term bounded above for all $r\in\widehat{S}^{(N),c}=[0,4N]$ as 
\begin{align*}
    |\textsc{Error}_{s,t}^{(N),d,c}|\leq CN^{-\chi}(1+\sqrt{r})^{\eta}.
\end{align*}
\end{proof}

\section{Convergence of the Laplace Transform}\label{s:p1}

In this section, we put the results from \cref{s:conv} together to prove \cref{l:laplaceconverge}. We begin by defining functions 
\begin{equation*}
    \begin{aligned}
        \mathcal{G}^{(N)}(\vec{r};\vec{c};\vec{x}) & :=\mathbf{1}_{\vec{r} \in \mathbb{R}^{d+1}_{\leq 4N}}\cdot 2^{-N}\prod_{k=1}^{d+1}\left(\cosh{\left(\frac{s_{k}}{\sqrt{N}}\right)}+1-\frac{r_{k}}{2N}\right)^{N(x_{k}^{(N)}-x_{k-1}^{(N)})}, \\ 
        \mathcal{G}(\vec{r};\vec{c};\vec{x}) & := \exp{\left(\frac{1}{4}\sum_{k=1}^{d}(s_{k}^{2}-r_{k})(x_{k}-x_{k-1})\right)}.
    \end{aligned}
\end{equation*} We make these definitions in order to concisely rewrite the formula for the Laplace transform of the diffusively scaled height function \eqref{e:laplace} in terms of $\mathcal{G}^{(N)}(\vec{r};\vec{c};\vec{x})$ and $\mathcal{G}(\vec{r};\vec{c};\vec{x})$,

\begin{align*}
    \varphi^{(N)}_{u,v,w,r}(\vec{c},\vec{x}) & :=\frac{\mathbb{E}\left[\prod_{k=1}^{d+1}\left(\cosh{\left(s_{k}^{(N)}\right)}+\mathbb{Y}_{e^{-2s_{k}^{(N)}}}\right)^{n_{k}-n_{k-1}}\right]}{\mathbb{E}\left[(1+\mathbb{Y}_{1})^{N}\right]} \\ & = \frac{\mathbb{E}\left[\mathcal{G}^{(N)}\left((\widehat{\mathbb{Y}}_{s_{1}}^{(N)},...,\widehat{\mathbb{Y}}^{(N)}_{s_{d+1}});\vec{c};\vec{x}\right)\right]}{\mathbb{E}\left[\mathcal{G}^{(N)}\left((\widehat{\mathbb{Y}}_{0}^{(N)},...,\widehat{\mathbb{Y}}^{(N)}_{0});\vec{0};\vec{x}\right)\right]}  =\frac{\mathbb{E}\left[\mathcal{G}^{(N)}\left((\widehat{\mathbb{Y}}_{s_{1}}^{(N)},...,\widehat{\mathbb{Y}}^{(N)}_{s_{d+1}});c;X\right)\right]}{2^{-N}\mathbb{E}\left[\left(1-\frac{\widehat{\mathbb{Y}}_{0}^{(N)}}{4N}\right)^{N}\right]} , \\ 
    \varphi_{u,v}(\vec{c};\vec{x}) & := \frac{\mathbb{E}\left[\exp{\left(\frac{1}{4}\sum_{k=1}^{d}\left(s_{k}^{2}-\mathbb{T}_{s_{k}}\right)\left(x_{k}-x_{k-1}\right)\right)}\right]}{\mathbb{E}\left[\left(-\frac{1}{4}\mathbb{T}_{0}\right)^{N}\right]} \\ & = \frac{\mathbb{E}\left[\mathcal{G}\left((\mathbb{T}_{s_{1}},...,\mathbb{T}_{s_{d}});\vec{c};\vec{x}\right)\right]}{\mathbb{E}\left[\mathcal{G}\left((\mathbb{T}_{0},...,\mathbb{T}_{0});\vec{c};\vec{x}\right)\right]}.
\end{align*}

In particular, we break this ratio up into its numerator and denominator, defining two associated functions

\begin{align*}
    \widetilde{\varphi}_{u,v,w,r}^{(N)}(\vec{c};\vec{x})& :=\int \mathcal{G}^{(N)}(\vec{r};\vec{c};\vec{x})\prod_{k=1}^{d}\widehat{\pi}^{(N)}_{s_{i+1},s_{i}}(r_{i+1},dr_{i})\cdot N^{u+v}\widehat{\pi}^{(N)}_{s_{d+1}}(dr_{d+1}) \\ & = N^{u+v}\mathbb{E}\left[\mathcal{G}^{(N)}\left((\widehat{\mathbb{Y}}_{s_{1}}^{(N)},...,\widehat{\mathbb{Y}}^{(N)}_{s_{d+1}});\vec{c};\vec{x}\right)\right]\\ 
    \widetilde{\varphi}_{u,v}(\vec{c};\vec{x}) &:= \int \mathcal{G}(\vec{r};\vec{c};\vec{x})\prod_{k=1}^{d}\mathfrak{p}_{s_{i+1},s_{i}}(r_{i+1},dr_{i})\cdot \mathfrak{p}_{s_{d+1}}(dr_{d+1}) 
    \\ & = \mathbb{E}\left[\mathcal{G}\left((\mathbb{T}_{s_{1}},...,\mathbb{T}_{s_{d}});\vec{c};\vec{x}\right)\right].
\end{align*}

In this language, we may write $$\varphi^{(N)}_{u,v,w,r}(\vec{c};\vec{x}) = \frac{\widetilde{\varphi}^{(N)}_{u,v,w,r}(\vec{c};\vec{x})}{\widetilde{\varphi}^{(N)}_{u,v,w,r}(\vec{0};\vec{x})}.$$

The proof of \cref{t:1} will proceed by demonstrating that (\cref{l:laconv}) $$\lim_{N\to\infty}\widetilde{\varphi}_{u,v,w,r}^{(N)}(\vec{c};\vec{x})=C_{u,v,w,r}\widetilde{\varphi}_{u,v}(\vec{c};\vec{x}).$$

We cite the following lemma from \cite{CK21}. 

\begin{lemma}[Lemma 8.8~\cite{CK21}] For all $N\in\mathbb{N}$, the functions $\widetilde{\varphi}_{u,v}(\vec{c},\vec{x})$, $\vec{c}\in (0,C_{u,r})^{d}$, and $\widetilde{\varphi}_{u,v}(\vec{0},\vec{x})$ take values in $(0,\infty)$.   
\end{lemma} The proof is identical to that in \cite{CK21}, so we only paraphrase the argument here. The first step is to apply the bound $\mathcal{G}(\vec{r};\vec{c};\vec{x})\leq Ce^{-(1-x_{d})r_{d+1}/4}$ and then repeatedly integrate over the variables $r_{i}$ to obtain $$\widetilde{\varphi}_{v,w,u,r}(\vec{c};\vec{x})\leq C\int_{\mathbb{R}}e^{-(1-x_{d})r}\mathfrak{p}_{0}(dr).$$ Then, it is necessary to apply further bounds to the asymptotic behavior of $\mathfrak{p}_{0}^{c}$ to show that this integral is finite. Showing that $\widetilde{\varphi}_{u,v}(\vec{c};\vec{x})$ is positive involves noting that $\mathcal{G}(\vec{r};\vec{c};\vec{x})\geq \mathcal{G}(\vec{2};\vec{c};\vec{x})$ for all $\vec{r}\in[1,2]^{d+1}$ and then using this bound to obtain the lower bound
\begin{align*}
    \widetilde{\varphi}_{u,v}(\vec{c};\vec{x})\geq \mathcal{G}(\vec{2};\vec{c};\vec{x})\int_{[1,2]^{d+1}}\prod_{i=1}^{d}\mathfrak{p}_{s_{i+1},s_{i}}^{c,c}(r_{i+1},r_{i})dr_{i}\cdot \mathfrak{p}_{s_{d+1}}^{c}(r_{d+1})dr_{d+1}.
\end{align*} and then completing the bound by using the explicit formulas for the marginal and transition measures to integrate this expression.

We will also use the following lemma, analogous to Lemma 8.2 in \cite{CK21}.
\begin{lemma}[Lemma 8.2~\cite{CK21}]\label{l:gconv}
    For every compact interval $I\subset \mathbb{R}$, there exists a constant $C>0$ such that for all $\vec{c}\in\mathbb{R}^{d+1}$,
    $$\mathcal{G}^{(N)}(\vec{r};\vec{c};\vec{x})\leq C\prod_{k=1}^{d+1}e^{-(x_{k}-x_{k-1})\frac{r_{k}}{4} + \frac{|r_{k}|}{2N}},$$ where $\vec{x}$ and $\vec{c}$ are defined in \eqref{e:assumptions}. Furthermore, for all $\vec{r}^{(N)}\in\mathbb{R}^{d+1}$ with $\lim_{N\to\infty}\vec{r}^{(N)}=\vec{r},$ $$\lim_{N\to\infty}\mathcal{G}^{(N)}(\vec{r}^{(N)};\vec{c};\vec{x})=\mathcal{G}(\vec{r};\vec{c};\vec{x}).$$
\end{lemma} 
The proof of this lemma is also unchanged from the proof of the analogous lemma in \cite{CK21}. It follows from an involved sequence of applying elementary identities and approximations.

Finally, we will show the convergence of the sequence  $\{\widetilde{\varphi}^{(N)}_{u,v,w,r}(\vec{c};\vec{x})\}_{N\in\mathbb{N}}$.  
\begin{lemma}\label{l:laconv}
    For $(q,A,B,C,D)$ and $h_{u,v,w,r}^{(N)}(x)$ satisfying the scaling conditions defined in \cref{d:scalenew}, and for all $\vec{c}\in(0,C_{u,r})^{d}$ and $\vec{c}=\vec{0}$,  $$\lim_{N\to\infty}\widetilde{\varphi}_{u,v,w,r}^{(N)}(\vec{c};\vec{x})=C_{u,v,w,r}\widetilde{\varphi}_{u,v}(\vec{c};\vec{x}).$$
\end{lemma}
The argument is unchanged from that given in \cite{CK21}, albeit using the lemmas we established in \cref{s:conv} for general $w$ and $r$ and the scaling conditions in \cref{d:scalenew}. For this reason, we only provide a summary of the argument here. 

There are two main steps. The first is to show pointwise convergence of the integrands. This is accomplished by combining the results of \cref{l:gconv} with those of \cref{l:singlepoint}, \cref{l:transprob}, and \cref{l:ddconv}. Then we must show dominated convergence to get the full result. Using the error bounds in \cref{l:singlepoint}, \cref{l:transprob}, and \cref{l:ddconv}, and the fact that $e^{-cr}e^{CN^{-\chi}(1+\sqrt{r})^{\eta}}\leq C'e^{-c'r}$ for $C'>0$ sufficiently large and $c'>0$ sufficiently small, and the bound that 
\begin{align*}
    \mathcal{G}^{(N)}(\vec{r};\vec{c};\vec{x})\leq Ce^{-c(r_{1}+...+r_{d+1})},
\end{align*} which follows from \cref{l:gconv}, we obtain the dominating function
\begin{equation}\label{e:dominated}
    \begin{aligned} \bigg|\mathcal{G}^{(N)}(\vec{r}; \vec{c}; \vec{x})\prod_{i=1}^{d}\widehat{\pi}_{s_{i+1},s_{i}}^{(N),a_{i+1},a_{i}}(r_{i+1},r_{i})\cdot \widehat{\pi}_{s_{d+1}}^{(N),a_{d+1}}(r_{d+1}) \bigg| \\ \leq C e^{-c(r_{1}+...+r_{d+1})}\prod_{i=1}^{d}\mathfrak{p}_{s_{i+1},s_{i}}^{a_{i+1},a_{i}}\left(r_{i+1},r_{i}\right)\cdot\mathfrak{p}_{s_{d+1}}^{a_{d+1}}(r_{d+1}).
    \end{aligned}
\end{equation}
As in \cite{CK21}, we insert the multiplicative factor $1=\prod_{i=1}^{d+1}\left(\mathbf{1}_{r_{i}\geq 0}+\mathbf{1}_{r_{i}<0}\right)$ into the function $\mathcal{G}^{(N)}(\vec r;\vec c;\vec x)$ in $\widetilde{\varphi}^{(N)}_{u,v,w,r}(\vec c;\vec x)$ and expand to obtain $2^{d+1}$ terms, corresponding to a choice of the atomic or continuous support at each $r_{i}$. Thus,
\begin{align*}
    \widetilde{\varphi}^{(N)}_{u,v,w,r}(\vec c,\vec x) & = \sum_{I\subset[[1,d+1]]}\sum_{r_{J}}\int_{r_{J}}\mathcal{G}(\vec r;\vec c;\vec x)\prod_{i=1}^{d}\widehat{\pi}^{(N),a_{i+1},a_{i}}_{s_{i+1},s_{i}}(r_{i+1},r_{i})\cdot N^{u+v}\widehat{\pi}^{(N),a_{d+1}}_{s_{d+1}}(r_{d+1}),
\end{align*} where $J$ is in the complement of $I$ in $[[1,d+1]]$ and the sum over $r_{I}$ is an $|I|$-fold sum over $r_{i}\in\widehat{S}^{(N),d}_{s_{i}}$ for $i\in I$, and the integral over $r_{J}$ is a $|J|$-fold integral of the $r$ variables with indices in $J$ over $\widehat{S}^{(N),c}$. We must show that for any choice of $I$, the corresponding sum over $r_{I}$ and integral over $r_{J}$ converge. We apply \cref{l:singlepoint,l:transprob,l:ddconv,l:atoms1} and the pointwise convergence of $\mathcal{G}^{(N)}(\vec r^{(N)};\vec c;\vec x)$ in \cref{l:gconv} to obtain pointwise convergence of each term. To obtain dominated convergence, we need to take the dominating function from \eqref{e:dominated} and show that, for all $I\subset [[1,d+1]]$, 
\begin{align*}
    \sum_{r_{I}}\int_{r_{J}}C'e^{-c'(r_{1}+...+r_{d+1})}\prod_{i=1}^{d}\mathfrak{p}_{s_{i+1},s_{i}}^{a_{i+1},a_{i}}(r_{i+1},r_{i})\cdot\mathfrak{p}_{s_{d+1}}^{a_{d+1}}(r_{d+1})<\infty,
\end{align*} which is accomplished by casework on the structure of the set $I$, see \cite[Section 8.4]{CK21}.

\section{Tightness and Proof of \cref{t:1}}\label{s:proof}

In this section, we discuss the tightness of the laws of $\{h_{u,v,w,r}^{(N)}(\cdot)-h_{u,v,w,r}^{(N)}(0)\}_{N\in\mathbb{N}}$, and then give a full proof of \cref{t:1}. A standard result in the study of stochastic processes (see, for instance, Karatzas and Shreve~\cite{KS}) tells us that if, for all $0\leq x,y\leq 1$ there exist $\alpha,\beta,\nu,C\in\mathbb{R}_{>0}$ such that 
\begin{align}\label{e:tightness}
\sup_{N\in\mathbb{N}}\mathbb{E}\left[|h_{u,v,w,r}^{(N)}(x)-h_{u,v,w,r}^{(N)}(y)|^{\alpha}\right]\leq C|x-y|^{1+\beta},
\end{align} and such that the process, evaluated at $0$, has bounded supremum (which is immediate for our process since $h^{(N)}_{u,v,w,r}(0)-h^{(N)}_{u,v,w,r}(0)=0$), then the collection of laws of $\{h_{u,v,w,r}^{(N)}(\cdot)-h_{u,v,w,r}^{(N)}(0)\}_{N\in\mathbb{N}}$ is tight on $(D[0,1],\mathcal{B}(D[0,1]))$. Since $h_{u,v,w,r}^{(N)}(0)-h_{u,v,w,r}^{(N)}(0)=0$ for all $N\in\mathbb{N}$, the proof of tightness reduces to demonstrating the second bound in \eqref{e:tightness}.

We will not prove the following lemma, as the argument is unchanged from that presented in \cite{CK21}. The result relies on an attractive coupling of the open ASEP, which is also stated and proved in \cite{CK21} Lemma 4.1.  The argument proceeds by sandwiching the measure of $h^{(N)}_{u,v,w,r}(x)-h^{(N)}_{u,v,w,r}(y)$ between two stationary measures, $h_{u,-u,w,r}^{(N)}(x)-h_{u,-u,w,r}^{(N)}(0)$ and $h_{-v,v,w,r}^{(N)}(x)-h^{(N)}_{-v,v,w,r}(0)$ which are known (see \cite{BW19}, Remark 2.4) to be in product Bernoulli form. 

\begin{lemma}[Proposition 4.2 (4)~\cite{CK21}]\label{l:tightcond} Consider the invariant measure for the open ASEP with parameters $A,B,C,D$ given by \eqref{e:params}. Denote the associated diffusively scaled height function increment process by $\{h^{(N)}_{u,v,w,r}(x)-h^{(N)}_{u,v,w,r}(0)\}_{x\in [0,1]}$, as before. Then for all $u,v,w,r\in\mathbb{R}$, we have the following H\"older bound. For all $\theta\in(0,1/2)$ and every $n\in\mathbb{Z}_{\geq 1}$ there exists a constant $C(\theta,n,u,v)>0$ such that for every $x,y\in[0,1]$ and every $N\in\mathbb{Z}_{\geq 1}$ $$\left\| h_{u,v,w,r}^{(N)}(x)-h_{u,v,w,r}^{(N)}(y) \right\|_{n} \leq C(\theta,n,u,v)|x-y|^{\theta},$$ where $\|\cdot\|_{n}:=\mathbb{E}[(\cdot)^{n}]^{1/n}$.
\end{lemma}
This completes the proof of tightness for \cref{t:1}. It is important to note that, while \cite{CK21} states \cref{l:tightcond} for the scaling conditions in \cref{d:scale}, the proof relies on an attractive coupling of ASEP and a sandwiching argument for the diffusively scaled random walks $h_{u,-v,w,r}^{(N)}(x)$ and $h_{-v,v,w,r}^{(N)}(x)$ which is unchanged when we allow $w$ and $r$ to vary. 

Finally, we can prove that our sequence of measures converges weakly to $H_{u,v}(\cdot)-H_{u,v}(0)$, completing the proof of \cref{t:1}. 
\begin{proof}[Proof of  \cref{t:1}]
    We must show that, when $u+v>0$ and $w,r>0$, the finite dimensional distributions of $\{h_{u,v,w,r}^{(N)}(\cdot)-h_{u,v,w,r}^{(N)}(0)\}$ converge in distribution to those of $H_{u,v}(\cdot)-H_{u,v}(0)$. Using a standard result in the convergence of continuous stochastic processes (see, for example, Karatzas and Shreve~\cite{KS}, Ch. 2 Theorem 4.15), convergence in finite dimensional distributions with tightness of the sequence, from \cref{t:1} $(1)$, implies that $\{h_{u,v,w,r}^{(N)}(\cdot)-h_{u,v,w,r}^{(N)}(0)\}\implies H_{u,v}(\cdot)-H_{u,v}(0)$ on $(D[0,1],\mathcal{B}(D[0,1])).$

    For any $d\in\mathbb{Z}_{\geq 1}$ and $0< x_{1}<...<x_{d}\leq 1$ and $x^{(N)} = N^{-1}\lfloor Nx\rfloor$, we use $\mathbb{P}^{(N)}$ to denote the law of the vector $\left(h_{u,v,w,r}^{(N)}(x_{1}^{(N)})-h_{u,v,w,r}^{(N)}(0), ..., h_{u,v,w,r}^{(N)}(x_{d}^{(N)})-h_{u,v,w,r}^{(N)}(0)\right)$. Consequently, using \cref{t:1} $(1)$, we know that every subsequence of this collection has a convergent sub-subsequence. We know that the corresponding sequence of Laplace transforms converges to a unique limit, which is the Laplace transform of the open KPZ stationary measure on an open interval (\cref{l:laplaceconverge}). Therefore, applying \cref{t:laplace}, we see that when $u+v>0$ and $w,r>0$, $\{h^{(N)}_{u,v,w,r}(\cdot)-h^{(N)}_{u,v,w,r}(0)\} \implies H_{u,v}(\cdot) - H_{u,v}(0).$
\end{proof}

\bibliographystyle{alpha}
\bibliography{openasep}

\end{document}